\documentclass[12pt]{amsart}

\usepackage{amsthm}

\usepackage{mathrsfs,mathtools,latexsym,eucal}
\usepackage{amscd,amsfonts,amssymb,amsmath,amsthm}
\usepackage{mathrsfs}
\usepackage{graphicx,graphics,color}
\usepackage[utf8]{inputenc}
\usepackage[active]{srcltx}
\usepackage[english]{babel}
\usepackage[pagewise,displaymath,mathlines]{lineno}
\usepackage{listings}
\usepackage{color}
\usepackage{subfigure} 
\usepackage{url}
\usepackage{float}
\usepackage{enumerate}

\usepackage[active]{srcltx}

\usepackage{pgf,tikz}
\usetikzlibrary{arrows}

\definecolor{dkgreen}{rgb}{0,0.6,0}
\definecolor{gray}{rgb}{0.5,0.5,0.5}
\definecolor{mauve}{rgb}{0.58,0,0.82}

\newtheorem*{theorem*}{Theorem}
\newtheorem{theorem}{Theorem}[section]
\newtheorem{lemma}[theorem]{Lemma}
\newtheorem{proposition}[theorem]{Proposition}
\newtheorem{corollary}[theorem]{Corollary}
\newtheorem{scholium}[theorem]{Scholium}

\newtheorem{example}[theorem]{Example}

\newcommand{\K}{\mathbb{K}}

\newcommand{\C}{\mathbb{C}}
\newcommand{\R}{\mathbb{R}}
\newcommand{\N}{\mathbb{N}}

\renewcommand{\leq}{\leqslant}
\renewcommand{\geq}{\geqslant}

\numberwithin{equation}{section}

\newcommand{\Bo}{\mathcal{B}}

\newcommand{\Un}{\mathcal{U}}
\newcommand{\ob}{\mathrm{ob}}

\newcommand{\suppv}{\mathrm{suppv}}

\newcommand{\B}{\mathsf{B}}

\newcommand{\E}{\mathsf{S}}

\newcommand{\argmax}{\mathrm{arg}\max}
\newcommand{\argmin}{\mathrm{arg}\min}

\newcommand{\ld}{\ell\mathit{d}}
\newcommand{\rd}{\mathit{rd}}

\begin{document}

\title[Exact solutions to a maxmin problem]{Exact solutions to the maxmin problem $\left\{\begin{array}{l} \max \|Ax\| \\  \|Bx\| \leq 1\end{array}\right.$}

\author{Soledad Moreno-Pulido$^0$}
\address{Department of Mathematics, College of Engineering, University of Cadiz, Puerto Real 11510, Spain (EU)}
\email{{\tt soledad.moreno@uca.es}}

\author{Francisco Javier Garcia-Pacheco$^0$}
\address{Department of Mathematics, College of Engineering, University of Cadiz, Puerto Real 11510, Spain (EU)}
\email{{\tt garcia.pacheco@uca.es}}

\author{Clemente Cobos-Sanchez}
\address{Department of Electronics, College of Engineering, University of Cadiz, Puerto Real 11510, Spain (EU)}
\email{{\tt clemente.cobos@uca.es}}

\author{Alberto Sanchez-Alzola}\thanks{corresponding author}
\address{Department of Statistics and Operation Research, College of Engineering, University of Cadiz, Puerto Real 11510, Spain (EU)}
\email[Corresponding Author]{{\tt alberto.sanchez@uca.es}\vspace{0.5cm} \newline
}


\date{}

\footnotetext{The authors have been supported by the Research Grant PGC-101514-B-100 awarded by the Spanish Ministry of Science, Innovation and Universities and partially funded by FEDER.}
\thanks{All authors contributed equally to the work}

\begin{abstract}
In this manuscript we provide an exact solution to the maxmin problem $\left\{\begin{array}{l} \max \|Ax\| \\ \|Bx\|\leq 1 \end{array}\right.$, where $A$ and $B$ are real matrices. This problem comes from a remodeling of $\left\{\begin{array}{l} \max \|Ax\| \\ \min \|Bx\| \end{array}\right.$, because the latter problem has no solution. Our mathematical method comes from the Abstract Operator Theory, whose strong machinery allows us to reduce the first problem to $\left\{\begin{array}{l} \max \|Cx\| \\ \|x\|\leq 1 \end{array}\right.$, which can be solved exactly by relying on supporting vectors. Finally, as appendices, we provide two applications of our solution: first, we construct a truly optimal minimum stored-energy Transcranian Magnetic Stimulation (TMS) coil, and second, we find an optimal geolocation involving statistical variables.
\end{abstract}

\maketitle

\section{Introduction}

\subsection{Scope}

Many problems in different disciplines like Physics, Statistics, Economics or Engineering can be modeled by using matrices and their norms (see for instance \cite{H,Y}). Here in this article we make use of supporting vectors to reformulate and solve problems and situations that commonly appear in the previously mentioned disciplines.

Supporting vectors are widely known in the literature of Geometry of Banach Spaces and Operator Theory. They are commonly known as the unit vectors at which an operator attains its norm. In the matrix setting, the supporting vectors of a matrix $A$ are the solutions of $$\displaystyle{\max_{\|x\|_2=1} \|Ax\|_2^2}.$$ Supporting vectors are topologically and geometrically studied in \cite{CSGPMPSM,GPNG}. In addition, generalized supporting vectors are defined and studied in \cite{CSGPMPSM,CSGPMPSA}. Again in the matrix setting, the generalized supporting vectors of a sequence of matrices $(A_i)_{i\in\N}$ are the solutions of $$\displaystyle{\max_{\|x\|_2=1}\sum_{i=1}^\infty \|A_ix\|_2^2}.$$ This optimization problem clearly generalizes the previous one.

A first application of supporting vectors was given in \cite{CSGPGRH} where a TMS coil was truly optimally designed. In that paper a three-component problem is stated but only the case of one component was solved. In \cite{CSGPMPSA} the three-component case was stated and solved by means of the generalized supporting vectors. Moreover, an optimal location problem using Principal Component Analysis is solved by means of generalized supporting vectors.

For other perspective on supporting vectors and generalized supporting vectors, we refer the reader to \cite{GPNG}.

\subsection{Novelties} In this subsection we intend to enumerate the novelties provided by this work:
\begin{enumerate}
\item We provide an exact solution of an optimization problem, not an heuristic method for approaching it. Specifically, we solve the maxmin problem \begin{equation}\label{maxmint}\left\{\begin{array}{l} \max \|Ax\| \\  \|Bx\|\leq 1 \end{array}\right.\end{equation}
\item A MATLAB code is provided for computing the solution to the maxmin problem.
\item Our solution applies to design truly optimal minimum stored-energy TMS coils and to find optimal geolocations involving statistical variables.
\item This is an interdisciplinary work that englobes pure abstract nontrivial theorems with their proofs and programming codes with their results to directly apply them to real-life situations.
\end{enumerate}

\subsection{Preliminaries}

A multiobjective optimization problem has the form $$P:=\left\{\begin{array}{l} \min f_j(x),\; 1\leq j \leq l \\ g_i(x) \leq b_i,\; 1\leq i\leq k\end{array}\right.$$ where $f_1,\dots,f_l,g_1,\dots,g_k: X\to \mathbb{R}$ are functions defined on a nonempty set $X$. Two special sets are associated to $P$, the feasible solutions of $P$ $$\mathrm{fea}(P):=\{x\in X: g_i(x) \leq b_i\; \forall\ 1\leq i\leq k\}$$ and the set of optimal solutions of $P$ $$\mathrm{sol}(P):=\{x\in \mathrm{fea}(P): f_j(x)\leq f_j(y)\; \forall y\in \mathrm{fea}(P)\;\forall\ 1\leq j\leq l\}.$$

Any multiobjective optimization problem can be rewritten as the intersection of optimization problems, that is, if $$P_j:=\left\{\begin{array}{l} \min f_j(x) \\ g_i(x) \leq b_i,\; 1\leq i\leq k\end{array}\right.$$ for $1\leq j\leq l$, then
\begin{itemize}
\item $P=P_1\wedge \cdots \wedge P_l$,
\item $\mathrm{fea}(P)=\mathrm{fea}(P_j)$ for all $1\leq j\leq l$, and 
\item $\mathrm{sol}(P)=\mathrm{sol}(P_1)\cap \cdots \cap \mathrm{sol}(P_l)$.
\end{itemize}
It commonly happens with multiobjective optimization problems that $$\mathrm{sol}(P)=\mathrm{sol}(P_1)\cap \cdots \cap \mathrm{sol}(P_l)=\varnothing.$$ In this situation, we have to search for another multiobjective optimization problem which has a solution and still models accurately the real-life situation where Problem $P$ comes from. In order to avoid the lack of solutions, it is a common practice to reduce the multiobjective optimization problem into a single optimization problem (increasing the number of constraints). Two typical reformulation are the following: 
\begin{equation}\label{refor1}
\left\{\begin{array}{l} \min f_j(x),\; 1\leq j \leq l \\ g_i(x) \leq b_i,\; 1\leq i\leq k\end{array}\right.\stackrel{\text{reform}}{\longrightarrow}  \left\{\begin{array}{l} \min f_{j_0}(x)\\ f_j(x)\leq c_j,\; 1\leq j \leq l, \; j\neq j_0\\ g_i(x) \leq b_i,\; 1\leq i\leq k\end{array}\right.
\end{equation}
and
\begin{equation}\label{refor2}
\left\{\begin{array}{l} \min f_j(x),\; 1\leq j \leq l \\ g_i(x) \leq b_i,\; 1\leq i\leq k\end{array}\right.\stackrel{\text{reform}}{\longrightarrow}  \left\{\begin{array}{l} \min h\left(f_1(x),\dots,f_l(x)\right)\\  g_i(x) \leq b_i,\; 1\leq i\leq k\end{array}\right.
\end{equation}
where $h:\R^l\to\R$ is a function conveniently chosen (usually an increasing function on each coordinate).

On the other hand, observe that if $\phi:Y\to X$ is a bijection, then it is easy to check that $\mathrm{fea}(P)=\phi(\mathrm{fea}(Q))$ and $\mathrm{sol}(P)=\phi(\mathrm{sol}(Q))$ where
\begin{equation}\label{cv1}
Q:=\left\{\begin{array}{l} \min (f_j\circ\phi)(y),\; 1\leq j \leq l \\ (g_i\circ \phi)(y) \leq b_i,\; 1\leq i\leq k\end{array}\right.
\end{equation}

Also note that $\mathrm{fea}(P)=\mathrm{fea}(R)$ and $\mathrm{sol}(P)=\mathrm{sol}(R)$ where
\begin{equation}\label{cv2}
R:=\left\{\begin{array}{l} \min (\phi_j\circ f_j)(x),\; 1\leq j \leq l \\ (\chi_i\circ g_i)(x) \leq \chi_i(b_i),\; 1\leq i\leq k\end{array}\right.
\end{equation}
and $\phi_j,\chi_i: \mathbb{R}\to \mathbb{R}$ are strictly increasing for $1\leq j\leq l$ and $1\leq i\leq k$.

The original maxmin optimization problem has the form $$M:=\left\{\begin{array}{l} \max g(x) \\ \min f(x) \end{array}\right.$$ where $f,g:X\to(0,\infty)$ are functions and $X$ is a nonempty set. Notice that $$\mathrm{sol}(M)=\argmax g(x) \cap \argmin f(x).$$ Many real-life problems can be mathematically model like a maxmin. However, this kind of multiobjective optimization problems may have the inconvenience of lacking a solution. If this occurs, then we are in need of remodeling the real-life problem with another mathematical optimization problem that has a solution and still models the real-life problem very accurately.

In \cite[Theorem 5.1]{CSGPGRH} it was shown that $$\argmax g(x) \cap \argmin f(x) \subseteq \argmin \frac{f(x)}{g(x)}\subseteq \{x\in X: \forall y\in X\; f(x)\leq f(y) \lor g(x)\geq g(y)\}.$$ This suggests that, in case $\mathrm{sol}(M)=\varnothing$, the following optimization problems are good alternatives to keep modeling the real-life problem accurately:
\begin{itemize}
\item $\left\{\begin{array}{l} \max g(x) \\ \min f(x) \end{array}\right.\stackrel{\text{reform}}{\longrightarrow} \left\{\begin{array}{l}\min \frac{f(x)}{g(x)}\\ g(x)\neq 0\end{array}\right.$. Here we have used the second typical reformulation for $h(u,v)=\frac{v}{u}$ described in Equation \eqref{refor2}.
\item $\left\{\begin{array}{l} \max g(x) \\ \min f(x) \end{array}\right.\stackrel{\text{reform}}{\longrightarrow} \left\{\begin{array}{l}\max \frac{g(x)}{f(x)}\\ f(x)\neq 0\end{array}\right.$. Here we have used the second typical reformulation for $k(u,v)=\frac{u}{v}$ described in Equation \eqref{refor2}.
\item $\left\{\begin{array}{l} \max g(x) \\ \min f(x) \end{array}\right.\stackrel{\text{reform}}{\longrightarrow} \left\{\begin{array}{l} \max g(x)\\ f(x)\leq a\end{array}\right.$. Here we have used the first typical reformulation \eqref{refor1} where $a$ is an appropriately chosen constant.
\item $\left\{\begin{array}{l} \max g(x) \\ \min f(x) \end{array}\right.\stackrel{\text{reform}}{\longrightarrow} \left\{\begin{array}{l} \min f(x)\\ g(x)\geq b\end{array}\right.$. Here we have used the first typical reformulation \eqref{refor1} where $b$ is an appropriately chosen constant.
\end{itemize}

We will prove in the third section that all four previous reformulations are equivalent for the original maxmin $\left\{\begin{array}{l} \max \|Ax\| \\ \min \|Bx\| \end{array}\right.$. In the fourth section, we will solve the reformulation $\left\{\begin{array}{l} \max \|Ax\| \\  \|Bx\| \leq 1 \end{array}\right.$.

\section{Characterizations of operators with null kernel}

Kernels will play a fundamental role towards solving the general reformulated maxmin \eqref{grmn} as shown in the next section. This is why we first study the operators with null kernel.

Throughout this section, all monoid actions considered will be left, all rngs will be associative, all rings will be unitary rngs, all absolute semi-values and all semi-norms will be non-zero, all modules over rings will be unital, all normed spaces will be real or complex and all algebras will be unitary and complex.

Given a rng $R$ and an element $s\in R$, we will denote by $\ld (s)$ to the set of left divisors of $s$, that is, $$\ld(s):=\{r\in R: \exists\; t\in R\setminus \{0\} \text{ with }rt=s\}.$$ Similarly, $\rd (s)$ stands for the set of right divisors of $s$. If $R$ is a ring, then the set of its invertibles is usually denoted by $\Un(R)$. Notice that $\ld(1)$ ($\rd(1)$) is precisely the subset of elements of $R$ which are right-(left) invertible. As a consequence, $\Un(R)=\ld(1)\cap \rd(1)$. Observe also that $\ld(0)\cap \rd(1)=\varnothing=\rd(0)\cap \ld(1)$. In general we have that $\ld(0)\cap \ld(1)\neq \varnothing$ and $\rd(0)\cap \rd(1)\neq \varnothing$. Later on in Example \ref{rd0rd1} we will provide an example of a ring where $\rd(0)\cap\rd(1)\neq \varnothing$.

Recall that an element $p$ of a monoid is called involutive if $p^2=1$. Given a rng $R$, an involution is an additive, antimultiplicative, composition-involutive map $*:R\to R$. A $*$-rng is a rng endowed with an involution.

The categorical concept of monomorphism will be very present in this work. Recall that a morphism $f\in\hom_{\mathcal{C}}(A,B)$ between objects $A$ and $B$ in a category $\mathcal{C}$ is called a monomorphism provided that $f\circ g=f\circ h$ implies $g=h$ for all $C\in\ob(\mathcal{C})$ and all $g,h\in\hom_{\mathcal{C}}(C,A)$. It is not hard to check that if $f\in\hom_{\mathcal{C}}(A,B)$ and there exist $C_0\in\ob(\mathcal{C})$ and $g_0\in\hom_{\mathcal{C}}(B,C_0)$ such that $g_0\circ f$ is a monomorphism, then $f$ is also a monomorphism. In particular, if $f\in\hom_{\mathcal{C}}(A,B)$ is a section, that is, exists $g\in\hom_{\mathcal{C}}(B,A)$ such that $g\circ f= I_A$, then $f$ is a monomorphism. As a consequence, the elements of the $\hom_{\mathcal{C}}(A,A)$ that have a left inverse are monomorphisms. In some categories, the last condition suffices to characterize monomorphisms. This is the case, for instance, of the category of vector spaces over a division ring.

Recall that $\mathcal{CL}(X,Y)$ stands for the space of continuous linear operators from $X$ to $Y$.

\begin{proposition}\label{aux}
A continuous linear operator $T:X\to Y$ between locally convex Hausdorff topological vector spaces $X$ and $Y$ verifies that $\ker(T)\neq \{0\}$ if and only if exists $S\in\mathcal{CL}(Y,X)\setminus\{0\}$ with $T\circ S=0$. In particular, if $X=Y$, then $\ker(T)\neq \{0\}$ if and only if $T\in \ld (0)$ in $\mathcal{CL}(X)$.
\end{proposition}

\begin{proof}
Let $S\in\mathcal{CL}(Y,X)\setminus\{0\}$ such that $T\circ S=0$. Fix any $y\in Y\setminus \ker(S)$, then $S(y)\neq 0$ and $T(S(y))=0$ so $S(y)\in \ker(T)\setminus \{0\}$. Conversely, if $\ker(T)\neq \{0\}$, then fix $x_0\in\ker(T)\setminus\{0\}$ and $y^*_0\in Y^*\setminus\{0\}$ (the existence of $y^*$ is guaranteed by the Hahn-Banach Theorem on the Hausdorff locally convex topological vector space $Y$). Next, consider $$\begin{array}{rrcl} S: & Y &\to & X\\ & y& \mapsto & S(y):=y^*_0(y)x_0.\end{array}$$Notice that $S\in\mathcal{CL}(Y,X)\setminus\{0\}$ and $T\circ S=0$.
\end{proof}

\begin{theorem}\label{aux2}
Let $X$ and $Y$ be locally convex Hausdorff topological vector spaces and $T:X\to Y$ a continuous linear operator.
\begin{enumerate}
\item If $T$ is a section, then $\ker(T)=\{0\}$
\item If $X$ and $Y$ are Banach spaces, $T(X)$ is complemented in $Y$ and $\ker(T)=\{0\}$, then $T$ is a section.
\end{enumerate}
\end{theorem}

\begin{proof}
\mbox{}
\begin{enumerate}
\item Trivial since sections are monomorphisms.
\item Consider $T:X\to T(X)$. Since $T(X)$ is complemented in $Y$ we have that it is closed in $Y$, thus it is a Banach space. Therefore, the Open Mapping Theorem assures that $T:X\to T(X)$ is an isomorphism. Let $T^{-1}:T(X)\to X$ be the inverse of $T:X\to T(X)$. Now consider $P:Y\to Y$ to be a continuous linear projection such that $P(Y)=T(X)$. Finally, it suffices to define $S:=T^{-1}\circ P$ since $S\circ T= I_X$.
\end{enumerate}
\end{proof}

We will finalize this section with a trivial example of a matrix $A\in\R^{3\times 2}$ such that $A\in\rd(I)\cap \rd(0)$.

\begin{example}\label{rd0rd1}
Consider $$A=\left(\begin{array}{cc} 1 & 0 \\ 0 & 1\\ 0 & 0\end{array}\right).$$ It is not hard to check that $\ker(A)=\{(0,0)\}$ thus $A$ is left-invertible by Theorem \ref{aux2}(2) and so $A\in \rd(I)$. In fact, $$\left(\begin{array}{ccc} 1 & 0 & 0 \\ 0 & 1  & 0\end{array}\right) \left(\begin{array}{cc} 1 & 0 \\ 0 & 1\\ 0 & 0\end{array}\right) = \left(\begin{array}{cc} 1 & 0 \\ 0 & 1\end{array}\right).$$ Finally, $$\left(\begin{array}{ccc} 0 & 0 & 1 \\ 0 & 0  & 1\end{array}\right) \left(\begin{array}{cc} 1 & 0 \\ 0 & 1\\ 0 & 0\end{array}\right) = \left(\begin{array}{cc} 0 & 0 \\ 0 & 0\end{array}\right).$$ 
\end{example}

\section{Remodeling the original maxmin problem $\left\{\begin{array}{l} \max \|T(x)\| \\ \min \|S(x)\| \end{array}\right.$ }

\subsection{The original maxmin problem has no solutions}

This subsection begins with the following theorem:

\begin{theorem}
If $T,S:X\to Y$ are nonzero continuous linear operators between Banach spaces $X$ and $Y$, then the maxmin problem
\begin{eqnarray}\label{gmn}
\left\{\begin{array}{l} \max \|T(x)\| \\ \min \|S(x)\| \end{array}\right.
\end{eqnarray}
has trivially no solution.
\end{theorem}

\begin{proof}
Observe that $\argmin\|S(x)\|=\ker(S)$ and $\argmax\|T(x)\|=\varnothing$ because $T\neq \{0\}$. Then the set of solutions of Problem \eqref{gmn} is $$\argmin\|S(x)\|\cap \argmax\|T(x)\|=\ker(S)\cap \varnothing=\varnothing.$$
\end{proof}

As a consequence, Problem \eqref{gmn} must be reformulated or remodeled.

\subsection{Equivalent reformulations for the original maxmin problem}

Following the first typical reformulation, given in Equation \eqref{refor1}, we obtain 
\begin{eqnarray}\label{grmn}
\left\{\begin{array}{l} \max \|T(x)\| \\ \|S(x)\|\leq 1 \end{array}\right.
\end{eqnarray}
Note that $\displaystyle \argmax_{\|S(x)\|\leq 1}\|T(x)\|$ is a $\mathbb{K}$-symmetric set, where $\mathbb{K}:=\R\text{ or }\C$, in other words, if $\lambda\in\K$ and $|\lambda|=1$, then $\lambda x\in \displaystyle{ \argmax_{\|S(x)\|\leq 1}\|T(x)\|}$ for every $x\in\displaystyle{ \argmax_{\|S(x)\|\leq 1}\|T(x)\|}$. The finite dimensional version of the previous reformulation is
\begin{eqnarray}\label{rmn}
\left\{\begin{array}{l} \max \|Ax\| \\ \|Bx\|\leq 1 \end{array}\right.
\end{eqnarray}
where $A,B\in\R^{m\times n}$.

Recall that $\Bo(X,Y)$ stands for the space of bounded operators from $X$ to $Y$.

\begin{lemma}\label{maxminchar1}
Let $X$ and $Y$ be Banach spaces and $T,S\in \Bo(X,Y)$. If the general reformulated maxmin problem $$\left\{\begin{array}{l} \max \|T(x)\| \\ \|S(x)\|\leq 1 \end{array}\right.$$ has a solution, then $\ker(S)\subseteq \ker(T)$.
\end{lemma}

\begin{proof}
If $\ker(S)\setminus \ker(T)\neq \varnothing$, then it suffices to consider the sequence $(nx_0)_{n\in\N}$ for $x_0\in\ker(S)\setminus\ker(T)$, since $\|S(nx_0)\|=0\leq 1$ for all $n\in \N$ and $\|T(nx_0)\|=n\|T(x_0)\|\to \infty$ as $n\to\infty$.
\end{proof}

The general maxmin \eqref{gmn} can also be reformulated by using the second typical reformulation \eqref{refor2}. This way we obtain $$\left\{\begin{array}{l} \max \|T(x)\| \\ \min \|S(x)\| \end{array}\right.\stackrel{\text{reform}}{\longrightarrow} \left\{\begin{array}{l} \max \frac{\|T(x)\|}{\|S(x)\|}\\ \|S(x)\|\neq 0 \end{array}\right.$$ 

\begin{lemma}\label{maxminchar2}
Let $X$ and $Y$ be Banach spaces and $T,S\in \Bo(X,Y)$. If the second general reformulated maxmin problem $$\left\{\begin{array}{l} \max \frac{\|T(x)\|}{\|S(x)\|}\\ \|S(x)\|\neq 0 \end{array}\right.$$ has a solution, then $\ker(S)\subseteq \ker(T)$.
\end{lemma}

\begin{proof}
Suppose there exists $x_0\in \ker(S)\setminus \ker(T)$. Then fix an arbitrary $x_1\in X\setminus \ker(S)$. Notice that $$\frac{\|T(nx_0+x_1)\|}{\|S(nx_0+x_1)\|}\geq\frac{n\|T(x_0)\|-\|T(x_1)\|}{\|S(x_1)\|}\to \infty$$ as $n\to \infty$.
\end{proof}

The next theorem shows that the previous two reformulations are in fact equivalent.

\begin{theorem}\label{maxminchar3}
Let $X$ and $Y$ be Banach spaces and $T,S\in \Bo(X,Y)$. Then $$\bigcup_{t>0}t\argmax_{\|S(x)\|\leq 1}\|T(x)\| = \argmax_{\|S(x)\|\neq 0} \frac{\|T(x)\|}{\|S(x)\|}.$$
\end{theorem}

\begin{proof}
Let $x_0\in \argmax_{\|S(x)\|\leq 1}\|T(x)\|$ and $t_0>0$. Fix an arbitrary $y\in X\setminus \ker(S)$. Notice that $x_0\notin \ker(S)$ in virtue of Theorem \ref{maxminchar1}. Then $$\|T(x_0)\|\geq \left\|T\left(\frac{y}{\|S(y)\|}\right)\right\|,$$ therefore $$\frac{\|T(tx_0)\|}{\|S(tx_0)\|}=\frac{\|T(x_0)\|}{\|S(x_0)\|}\geq \|T(x_0)\|\geq \left\|T\left(\frac{y}{\|S(y)\|}\right)\right\|.$$ Conversely, let $x_0\in \argmax_{\|S(x)\|\neq 0} \frac{\|T(x)\|}{\|S(x)\|}$. Fix an arbitrary $y\in X$ with $\|S(y)\|\leq 1$. Then $$\left\|T\left(\frac{x_0}{\|S(x_0)\|}\right)\right\|=\frac{\|T(x_0)\|}{\|S(x_0)\|}\geq \frac{\|T(y)\|}{\|S(y)\|}\geq \|T(y)\|$$ which means that $$\frac{x_0}{\|S(x_0)\|}\in \argmax_{\|S(x)\|\leq 1}\|T(x)\|$$ and thus $$x_0\in\|S(x_0)\|\argmax_{\|S(x)\|\leq 1}\|T(x)\|\subseteq \bigcup_{t>0}t\argmax_{\|S(x)\|\leq 1}\|T(x)\|.$$
\end{proof}

The reformulation $$\left\{\begin{array}{l}\min \frac{\|S(x)\|}{\|T(x)\|}\\ \|T(x)\|\neq 0\end{array}\right.$$ is slightly different from the previous two reformulations. In fact, if $\ker(S)\setminus \ker(T)\neq \varnothing$, then $\argmin_{\|T(x)\|\neq 0} \frac{\|S(x)\|}{\|T(x)\|}=\ker(S)\setminus\ker(T)$. The previous reformulation is equivalent to the following one as shown in the next theorem: $$\left\{\begin{array}{l}\min  \|S(x)\| \\ \|T(x)\|\geq 1\end{array}\right.$$

\begin{theorem}
Let $X$ and $Y$ be Banach spaces and $T,S\in \Bo(X,Y)$. Then $$\bigcup_{t>0}t\argmin_{\|T(x)\|\geq 1}\|S(x)\| = \argmin_{\|T(x)\|\neq 0} \frac{\|S(x)\|}{\|T(x)\|}.$$
\end{theorem}

We spare of the details of the proof of the previous theorem to the reader. Notice that if $\ker(S)\setminus \ker(T)\neq \varnothing$, then $\argmin_{\|T(x)\|\geq 1}  \|S(x)\| =\ker(S) \setminus \{x\in: \|T(x)\|< 1\}$. However, if $\ker(S)\subseteq \ker(T)$, then all four reformulations are equivalent, as shown in the next theorem, whose proof’s details we spare again to the reader.

\begin{theorem}\label{equa}
Let $X$ and $Y$ be Banach spaces and $T,S\in\Bo(X,Y)$. If $\ker(S)\subseteq \ker(T)$, then $$\argmax_{\|S(x)\|\neq 0}\frac{\|T(x)\|}{\|S(x)\|}=\argmin_{\|T(x)\|\neq 0}\frac{\|S(x)\|}{\|T(x)\|}.$$
\end{theorem}

\section{Solving the maxmin problem $\left\{\begin{array}{l} \max \|T(x)\| \\  \|S(x)\|\leq 1 \end{array}\right.$ }

We will distinguish between two cases.

\subsection{First case: $S$ is an isomorphism over its image}

By bearing in mind Theorem \ref{equa}, we can focus on the first reformulation proposed at the beginning of the previous section: $$\left\{\begin{array}{l} \max \|T(x)\| \\ \min \|S(x)\| \end{array}\right.\stackrel{\text{reform}}{\longrightarrow} \left\{\begin{array}{l} \max \|T(x)\|\\ \|S(x)\|\leq 1\end{array}\right.$$

The idea we propose to solve the previous reformulation is to make use of supporting vectors (see \cite{CSGPGRH,CSGPMPSM,CSGPMPSA,GPNG}). Recall that if $R:X\to Y$ is a continuous linear operator between Banach spaces, then the set of supporting vectors of $R$ is defined by $$\suppv(R):=\argmax_{\|x\|\leq 1}\|R(x)\|.$$ The idea of using supporting vectors is that the optimization problem $$\left\{\begin{array}{l} \max \|R(x)\| \\ \|x\|\leq 1 \end{array}\right.$$ whose solutions are by definition the supporting vectors of $R$, can be easily solved theoretically and computationally (see \cite{CSGPMPSA}).

Our first result towards this direction considers the case where $S$ is an isomorphism over its image.

\begin{theorem}\label{th1}
Let $X$ and $Y$ be Banach spaces and $T,S\in\Bo(X,Y)$. Suppose that $S$ is an isomorphism over its image and $S^{-1}:S(X)\to X$ denotes its inverse. Suppose also that $S(X)$ is complemented in $Y$, being $p:Y\to Y$ a continuous linear projection onto $S(X)$. Then $$S^{-1}\left(S(X)\cap \argmax_{\|y\|\leq 1} \left\|\left(T\circ S^{-1}\circ p\right)(y)\right\|\right)\subseteq \argmax_{\|S(x)\|\leq 1}\|T(x)\|.$$ If, in addition, $\|p\|=1$, then $$\argmax_{\|S(x)\|\leq 1}\|T(x)\|= S^{-1}\left(S(X)\cap \argmax_{\|y\|\leq 1} \left\|\left(T\circ S^{-1}\circ p\right)(y)\right\|\right).$$ 
\end{theorem}

\begin{proof}
We will show first that $$S(X)\cap \argmax_{\|y\|\leq 1} \left\|\left(T\circ S^{-1}\circ p\right)(y)\right\|\subseteq S\left(\argmax_{\|S(x)\|\leq 1}\|T(x)\|\right).$$ Let $y_0=S(x_0)\in \displaystyle{\argmax_{\|y\|\leq 1} \left\|\left(T\circ S^{-1}\circ p\right)(y)\right\|}$. We will show that $x_0\in \displaystyle{\argmax_{\|S(x)\|\leq 1}\|T(x)\|}$. Indeed, let $x\in X$ with $\|S(x)\|\leq 1$. Since $\|S(x_0)\|=\|y_0\|\leq 1$, by assumption we obtain
\begin{eqnarray*}
\|T(x)\|&=&\left\|\left(T\circ S^{-1}\circ p\right)(S(x))\right\|\\
&\leq & \left\|\left(T\circ S^{-1}\circ p\right)(y_0)\right\|\\
&=&\left\|\left(T\circ S^{-1}\circ p\right)(S(x_0))\right\|\\
&=&\left\|T(x_0)\right\|.
\end{eqnarray*}
Now assume that $\|p\|=1$. We will show that $$S\left(\argmax_{\|S(x)\|\leq 1}\|T(x)\|\right)\subseteq S(X)\cap \argmax_{\|y\|\leq 1} \left\|\left(T\circ S^{-1}\circ p\right)(y)\right\|.$$ Let $x_0\in \displaystyle{\argmax_{\|S(x)\|\leq 1}\|T(x)\|}$, we will show that $S(x_0)\in \displaystyle{\argmax_{\|y\|\leq 1} \left\|\left(T\circ S^{-1}\circ p\right)(y)\right\|}$. Indeed, let $y\in\B_Y$. Observe that $$\left\|S\left( S^{-1}( p(y))\right)\right\|=\|p(y)\|\leq \|y\|\leq 1$$ so by assumption
\begin{eqnarray*}
\left\|\left(T\circ S^{-1}\circ p\right)(y)\right\|&=&\left\|T\left(S^{-1}( p(y))\right)\right\|\\
&\leq& \|T(x_0)\|\\
&=&\left\|T\left(S^{-1}(p(S(x_0)))\right)\right\|\\
&=&\left\|\left(T\circ S^{-1}\circ p\right) (S(x_0))\right\|.
\end{eqnarray*}
\end{proof}

Notice that, in the settings of Theorem \ref{th1}, $S^{-1}\circ p$ is a left-inverse of $S$, in other words, $S$ is a section, as in Theorem \ref{aux2}(2).

Taking into consideration that every closed subspace of a Hilbert space is $1$-complemented (see \cite{Bo,Ka} to realize that this fact characterizes Hilbert spaces of dimension $\geq 3$), we directly obtain the following corollary.

\begin{corollary}\label{cor1}
Let $X$ be a Banach space, $Y$ a Hilbert space and $T,S\in\Bo(X,Y)$ such that $S$ is an isomorphism over its image and $S^{-1}:S(X)\to X$ its inverse. Then 
\begin{eqnarray*}
\argmax_{\|S(x)\|\leq 1}\|T(x)\|&=& S^{-1}\left(S(X)\cap \argmax_{\|y\|\leq 1} \left\|\left(T\circ S^{-1}\circ p\right)(y)\right\|\right)\\
&=& S^{-1}\left(S(X)\cap \suppv\left(T\circ S^{-1}\circ p\right)\right)
\end{eqnarray*}
where $p:Y\to Y$ is the orthogonal projection on $S(X)$.
\end{corollary}

\subsection{The Moore-Penrose inverse}

If $B\in \K^{m\times n}$, then the Moore-Penrose inverse of $B$, denoted by $B^+$, is the only matrix $B^+\in\K^{n\times m}$ which verifies the following:
\begin{itemize}
    \item $BB^+B=B$.
    \item $B^+BB^+=B^+$.
    \item $(BB^+)^*=BB^+$.
    \item $(B^+B)^*=B^+B$.
\end{itemize}

If $\ker(B)=0$, then $B^+$ is a left-inverse of $B$. Even more, $BB^+$ is the orthogonal projection onto the range of $B$, thus we have the following scholium from Corollary \ref{cor1}.

\begin{scholium}
Let $A,B\in\R^{m\times n}$ such that $\ker(B)=\{0\}$. Then
\begin{eqnarray*}
B\left(\argmax_{\|Bx\|_2\leq 1}\|Ax\|_2\right)&=& B\R^n \cap \argmax_{\|y\|_2\leq 1} \left\|AB^+y\right\|_2 \\
&=&  B\R^n \cap \suppv \left(AB^+\right)
\end{eqnarray*}
\end{scholium}

According to the previous scholium, in its settings, if $y_0\in \argmax_{\|y\|_2\leq 1} \left\|AB^+y\right\|_2$ and there exists $x_0\in \R^n$ such that $y_0=Bx_0$, then $x_0\in \argmax_{\|Bx\|_2\leq 1}\|Ax\|_2$ and $x_0$ can be computed as $$x_0=B^+Bx_0=B^+y_0.$$

\subsection{Second case: $S$ is not an isomorphism over its image}

What happens if $S$ is not an isomorphism over its image? Next theorem answers this question.

\begin{theorem}\label{th2}
Let $X$ and $Y$ be Banach spaces and $T,S\in\Bo(X,Y)$ such that $\ker(S)\subseteq \ker(T)$. If $$\begin{array}{rrcl}\pi:&X&\to& X/\ker(S)\\ &x&\mapsto&\pi(x):=x+\ker(S)\end{array}$$ denotes the quotient map, then $$\argmax_{\|S(x)\|\leq 1}\|T(x)\|=\pi^{-1}\left( \argmax_{\|\overline{S}(\pi(x))\|\leq 1}\|\overline{T}(\pi(x))\|\right),$$ where $$\begin{array}{rrcl}\overline{T}: & \frac{X}{\ker(S)}&\to&Y\\ & \pi(x)&\mapsto&\overline{T}(\pi(x)):=T(x)\end{array}$$ and $$\begin{array}{rrcl}\overline{S}: & \frac{X}{\ker(S)}&\to&Y\\ & \pi(x)&\mapsto&\overline{S}(\pi(x)):=S(x).\end{array}$$
\end{theorem}

\begin{proof}
Let $x_0\in \argmax_{\|S(x)\|\leq 1}\|T(x)\|$. Fix an arbitrary $y\in X$ with $\|\overline{S}(\pi(y))\|\leq 1$. Then $\|S(y)\|=\|\overline{S}(\pi(y))\|\leq 1$ therefore $$\|\overline{T}(\pi(x_0)\|=\|T(x_0)\|\geq \|T(y)\|=\|\overline{T}(\pi(y))\|.$$ This shows that $\pi(x_0)\in \argmax_{\|\overline{S}(\pi(x))\|\leq 1}\|\overline{T}(\pi(x))\|$. Conversely, let $$\pi(x_0)\in \argmax_{\|\overline{S}(\pi(x))\|\leq 1}\|\overline{T}(\pi(x))\|.$$ Fix an arbitrary $y\in X$ with $\|S(y)\|\leq 1$. Then $\|\overline{S}(\pi(y))\|=\|S(y)\|\leq 1$ therefore $$\|T(x_0)\|=\|\overline{T}(\pi(x_0))\|\geq \|\overline{T}(\pi(y))\|=\|T(y)\|.$$ This shows that $x_0\in\argmax_{\|S(x)\|\leq 1}\|T(x)\|$.
\end{proof}

Notice that, in the settings of Theorem \ref{th2}, if $S(X)$ is closed in $Y$, then $\overline{S}$ is an isomorphism over its image $S(X)$, and thus in this case Theorem \ref{th2} reduces the reformulated maxmin to Theorem \ref{th1}.

\subsection{Characterizing when the finite dimensional reformulated maxmin has a solution}

The final part of this section is aimed at characterizing when the finite dimensional reformulated maxmin has a solution.

\begin{lemma}\label{boundseq}
Let $S:X\to Y$ be a linear operator between finite dimensional Banach spaces $X$ and $Y$. If $(x_n)_{n\in \N}$ is a sequence in $\{x\in X: \|S(x)\|\leq 1\}$, then there exists a sequence $(z_n)_{n\in\N}$ in $\ker(S)$ such that $(x_n+z_n)_{n\in\N}$ is bounded.
\end{lemma}

\begin{proof}
Consider the linear operator
$$\begin{array}{rrcl}
\overline{S}:& \frac{X}{\ker(S)}&\to & Y\\
& x+ \ker(S) &\mapsto & \overline{S}(x+\ker(S))=S(x).
\end{array}$$
Note that $$\left\|\overline{S}(x_n+\ker(S))\right\|=\|S(x_n)\|\leq 1$$ for all $n\in\N$, therefore the sequence $(x_n+\ker(S))_{n\in\N}$ is bounded in $\frac{X}{\ker(S)}$ because $\frac{X}{\ker(S)}$ is finite dimensional and $\overline{S}$ has null kernel so its inverse is continuous. Finally, choose $z_n\in\ker(S)$ such that $\|x_n+z_n\|< \|x_n+\ker(S)\|+\frac{1}{n}$ for all $n\in\N$.
\end{proof}

\begin{lemma}\label{boundA}
Let $A,B\in \R^{m\times n}$. If $\ker(B)\subseteq \ker(A)$, then $A$ is bounded on $\{x\in \R^n: \|Bx\|\leq 1\}$ and attains its maximum on that set.
\end{lemma}

\begin{proof}
Let $(x_n)_{n\in \N}$ be a sequence in $\{x\in \R^n: \|Bx\|\leq 1\}$. In accordance with Lemma \ref{boundseq}, there exists a sequence $(z_n)_{n\in\N}$ in $\ker(B)$ such that $(x_n+z_n)_{n\in\N}$ is bounded. Since $A(x_n)=A(x_n+z_n)$ by hypothesis (recall that $\ker(B)\subseteq \ker(A)$), we conclude that $A$ is bounded on $\{x\in \R^n: \|Bx\|\leq 1\}$. Finally, let $(x_n)_{n\in\N}$ be a sequence in $\{x\in \R^n: \|Bx\|\leq 1\}$ such that $\|Ax_n\|\to \displaystyle{\max_{\|Bx\|\leq 1}\|Ax\|}$ as $n\to\infty$. Note that $\left\|\overline{A}(x_n+\ker(B))\right\|=\|Ax_n\|$ for all $n\in\N$, so $\left(\overline{A}(x_n+\ker(B))\right)_{n\in\N}$ is bounded in $\R^m$ and so is $\left(\overline{A}(x_n+\ker(B))\right)_{n\in\N}$ in $\frac{\R^n}{\ker(B)}$. Fix $b_n\in\ker(B)$ such that $\|x_n+b_n\|<\|x_n+\ker(B)\|+\frac{1}{n}$ for all $n\in \N$. This means that $(x_n+b_n)_{n\in\N}$ is a bounded sequence in $\R^n$ so we can extract a convergent subsequence $\left(x_{n_k}+b_{n_k}\right)_{k\in\N}$ to some $x_0\in X$. At this stage, notice that $\left\|B\left(x_{n_k}+b_{n_k}\right)\right\|=\left\|Bx_{n_k}\right\|\leq 1$ for all $k\in\N$ and $\left(B\left(x_{n_k}+b_{n_k}\right)\right)_{k\in\N}$ converges to $Bx_0$, so $\|Bx_0\|\leq 1$. Note also that, since $\ker(B)\subseteq \ker(A)$, $\left(\left\|Ax_{n_k}\right\|\right)_{n\in\N}$ converges to $\|Ax_0\|$, which implies that $$x_0\in \argmax_{\|Bx\|\leq 1}\|Ax\|.$$
\end{proof}

\begin{theorem}\label{maxminchar}
Let $A,B\in \R^{m\times n}$. The reformulated maxmin problem $$\left\{\begin{array}{l} \max \|Ax\| \\ \|Bx\|\leq 1 \end{array}\right.$$ has a solution if and only if $\ker(B)\subseteq \ker(A)$.
\end{theorem}

\begin{proof}
If $\ker(B)\subseteq \ker(A)$, then we just need to call on Lemma \ref{boundA}. Conversely, if $\ker(B)\setminus \ker(A)\neq \varnothing$, then it suffices to consider the sequence $(nx_0)_{n\in\N}$ for $x_0\in\ker(B)\setminus\ker(A)$, since $\|B(nx_0)\|=0\leq 1$ for all $n\in \N$ and $\|A(nx_0)\|=n\|A(x_0)\|\to \infty$ as $n\to\infty$.
\end{proof}

\subsection{Matrices on quotient spaces}

Consider the maxmin $$\left\{\begin{array}{l} \max \|T(x)\|\\ \|S(x)\|\leq 1\end{array}\right.$$ where $X$ and $Y$ are Banach spaces and $T,S\in\Bo(X,Y)$ with $\ker(S)\subseteq \ker(T)$. Notice that if $(e_i)_{i\in I}$ is a Hamel basis of $X$, then $\left(e_i+\ker(S)\right)_{i\in I}$ is a generator system of $\frac{X}{\ker(S)}$. By making use of the Zorn’s Lemma, it can be shown that $\left(e_i+\ker(S)\right)_{i\in I}$ contains a Hamel basis of $\frac{X}{\ker(S)}$. Observe that a subset $C$ of $\frac{X}{\ker(S)}$ is linearly independent if and only if $S(C)$ is a linearly independent subset of $Y$.

In the finite dimensional case, we have $$\begin{array}{rrcl}\overline{B}: & \frac{\R^n
}{\ker(B)}&\to&\R^m\\ & x+\ker(B)&\mapsto&\overline{B}(x+\ker(B)):=Bx.\end{array}$$ and $$\begin{array}{rrcl}\overline{A}: & \frac{\R^n
}{\ker(B)}&\to&\R^m\\ & x+\ker(B)&\mapsto&\overline{A}(x+\ker(B)):=Ax.\end{array}$$ If $\{e_1,\dots,e_n\}$ denotes the canonical basis of $\R^n$, then $\{e_1+\ker(B),\dots,e_n+\ker(B)\}$ is a generator system of $\frac{\R^n
}{\ker(B)}$. This generator system contains a basis of $\frac{\R^n
}{\ker(B)}$ so let $\{e_{j_1}+\ker(B),\dots,e_{j_l}+\ker(B)\}$ be a basis of $\frac{\R^n
}{\ker(B)}$. Note that $\overline{A}\left(e_{j_k}+\ker(B)\right)=Ae_{j_k}$ and $\overline{B}\left(e_{j_k}+\ker(B)\right)=Be_{j_k}$ for every $k\in\{1,\dots,l\}$. Therefore, the matrix associated to the linear map defined by $\overline{B}$ can be obtained from the matrix $B$ by removing the columns corresponding to the indices $\{1,\dots,n\}\setminus \{j_1,\dots,j_l\}$, in other words, the matrix associated to $\overline{B}$ is $\left[Be_{j_1}|\cdots|Be_{j_l}\right]$. Similarly, the matrix associated to the linear map defined by $\overline{A}$ is $\left[Ae_{j_1}|\cdots|Ae_{j_l}\right]$. As we mentioned above, recall that a subset $C$ of $\frac{\R^n}{\ker(B)}$ is linearly independent if and only if $B(C)$ is a linearly independent subset of $\R^m$. As a consequence, in order to obtain the basis $\{e_{j_1}+\ker(B),\dots, e_{j_l}+\ker(B)\}$, it suffices to look at the rank of $B$ and consider the columns of $B$ that allow such rank, which automatically gives us the matrix associated to $\overline{B}$, that is, $\left[Be_{j_1}|\cdots|Be_{j_l}\right]$.

Finally, let $$\begin{array}{rrcl} \pi:&\R^n&\to&\frac{\R^n}{\ker(B)}\\ &x&\mapsto&\pi(x):x+\ker(B)\end{array}$$ denote the quotient map. Let $l:=\mathrm{rank}(B)=\dim\left(\frac{\R^n}{\ker(B)}\right)$. If $x=(x_1,\dots,x_l)\in\R^l$, then $\sum_{k=1}^lx_k\left(e_{j_k}+\ker(B)\right)\in \frac{\R^n}{\ker(B)}$. The vector $z\in\R^n$ defined by $$z_p:=\left\{\begin{array}{rl} x_k&p=j_k\\ 0 & p\notin\{j_1,\dots,j_l\}\end{array}\right.$$ verifies that $$p(z)=\sum_{k=1}^lx_k\left(e_{j_k}+\ker(B)\right).$$ To simplify the notation, we can define the map $$\begin{array}{rrcl}\alpha:&\R^l&\to&\R^n\\ &x&\mapsto&\alpha(x):=z\end{array}$$ where $z$ is the vector described right above.

\subsection{Conclusions: schematic summary}\label{css}

This subsection compiles all the results from the previous subsections and defines the structure of the algorithm that solves the maxmin.

Let $A,B\in\R^{m\times n}$ with $\ker(B)\subseteq \ker(A)$. Then $$\left\{\begin{array}{l} \max \|Ax\|_2 \\ \min \|Bx\|_2 \end{array}\right.\stackrel{\text{reform}}{\longrightarrow} \left\{\begin{array}{l} \max \|Ax\|_2\\ \|Bx\|_2\leq 1\end{array}\right.$$

\begin{itemize}
    \item[Case 1:] $\ker(B)=\{0\}$. $B^+$ denotes the Moore-Penrose inverse of $B$. $$\left\{\begin{array}{l} \max \|Ax\|_2 \\ \|Bx\|_2\leq 1 \end{array}\right.\stackrel{\text{supp. vec.}}{\longrightarrow} \left\{\begin{array}{l} \max \|AB^+y\|_2\\ \|y\|_2\leq 1\end{array}\right.\stackrel{\text{solution}}{\longrightarrow}\left\{\begin{array}{l}\displaystyle{y_0\in\argmax_{\|y\|_2\leq 1}\|AB^+y\|_2}\\ \mathrm{rank}(B)=\mathrm{rank}([B|y_0])\end{array}\right.\stackrel{\text{final sol.}}{\longrightarrow}x_0:=B^+y_0$$
    
    \item[Case 2:] $\ker(B)\neq \{0\}$. $\overline{B}=\left[Be_{j_1}|\cdots|Be_{j_l}\right]$ where $\mathrm{rank}(B)=l=\mathrm{rank}\left(\overline{B}\right)$ and $\overline{A}=\left[Ae_{j_1}|\cdots|Ae_{j_l}\right]$. $$\left\{\begin{array}{l} \max \|Ax\|_2 \\ \|Bx\|_2\leq 1 \end{array}\right.\stackrel{\text{case 1}}{\longrightarrow}\left\{\begin{array}{l} \max \|\overline{A}y\|_2 \\ \|\overline{B}y\|_2\leq 1 \end{array}\right.\stackrel{\text{solution}}{\longrightarrow}  y_0\stackrel{\text{final sol.}}{\longrightarrow}x_0:=\alpha(y_0)$$
\end{itemize}

\section{Maxmin involving more operators}\label{moreope}

Let $X$ and $Y$ be Banach spaces and $(T_n)_{n\in\N}$ and $(S_n)_{n\in\N}$ sequences of continuous linear operators from $X$ to $Y$. The maxmin
\begin{equation}\label{supermaxmin}
\left\{\begin{array}{l} 
\max \|T_n(x)\|\; n\in\N\\
\min \|S_n(x)\|\; n\in\N
\end{array}\right.
\end{equation}
can be reformulated like (recall the second typical reformulation)
\begin{equation}\label{supermaxmin2}
\left\{\begin{array}{l} 
\max \sum_{n=1}^\infty\|T_n(x)\|^2\\
\min \sum_{n=1}^\infty\|S_n(x)\|^2
\end{array}\right.
\end{equation}
which can be transformed into a regular maxmin like in \eqref{gmn} by considering the operators $$\begin{array}{rrcl} T:&X&\to&\ell_2(Y)\\ &x&\mapsto&T(x):=(T_n(x))_{n\in\N}\end{array}$$ and $$\begin{array}{rrcl} S:&X&\to&\ell_2(Y)\\ &x&\mapsto&S(x):=(S_n(x))_{n\in\N}\end{array}$$ obtaining then $$\left\{\begin{array}{l} 
\max \|T(x)\|^2\\
\min \|S(x)\|^2
\end{array}\right.$$ which is equivalent to $$\left\{\begin{array}{l} 
\max \|T(x)\|\\
\min \|S(x)\|
\end{array}\right.$$
Observe that for the operators $T$ and $S$ to be well defined it is sufficient that $(\|T_n\|)_{n\in\N}$ and $(\|S_n\|)_{n\in\N}$ be in $\ell_2$.

\appendix

\section{Applications to optimal TMS coils}

\subsection{Introduction to TMS coils}

Transcranial Magnetic Stimulation (TMS) is a non-invasive technique to stimulate the brain, which is applied to  psychiatric and medical conditions, such as major depressive disorder, schizophrenia, bipolar depression, post-traumatic, stress disorder and obsessive-compulsive disorder, amongst others \cite{key:wassermann}. In TMS, strong current pulses driven through a coil are used to induce an electric field stimulating neurons in the cortex. 

The goal in TMS coil design is to find optimal positions for the multiple windings of coils (or equivalently the  current density) so as to produce fields with the desired spatial characteristics and properties \cite{CSGPGRH} (high focality, field penetration depth, low inductance, etc.), where this design problem has been frequently posed as a constrained optimization problem.

Moreover, an important safety issue in TMS is the minimization of the stimulation of non-target areas. Therefore, the development of TMS as a medical tool would be benefited with the design of TMS stimulators capable of inducing a maximum electric field in the region of interest, while minimizing the undesired stimulation in other prescribed regions.

\subsection{Minimum stored-energy TMS coil}

\begin{figure}[h]
\begin{center}
\subfigure[] {\label{fig:rois}\includegraphics[width=0.9\textwidth]{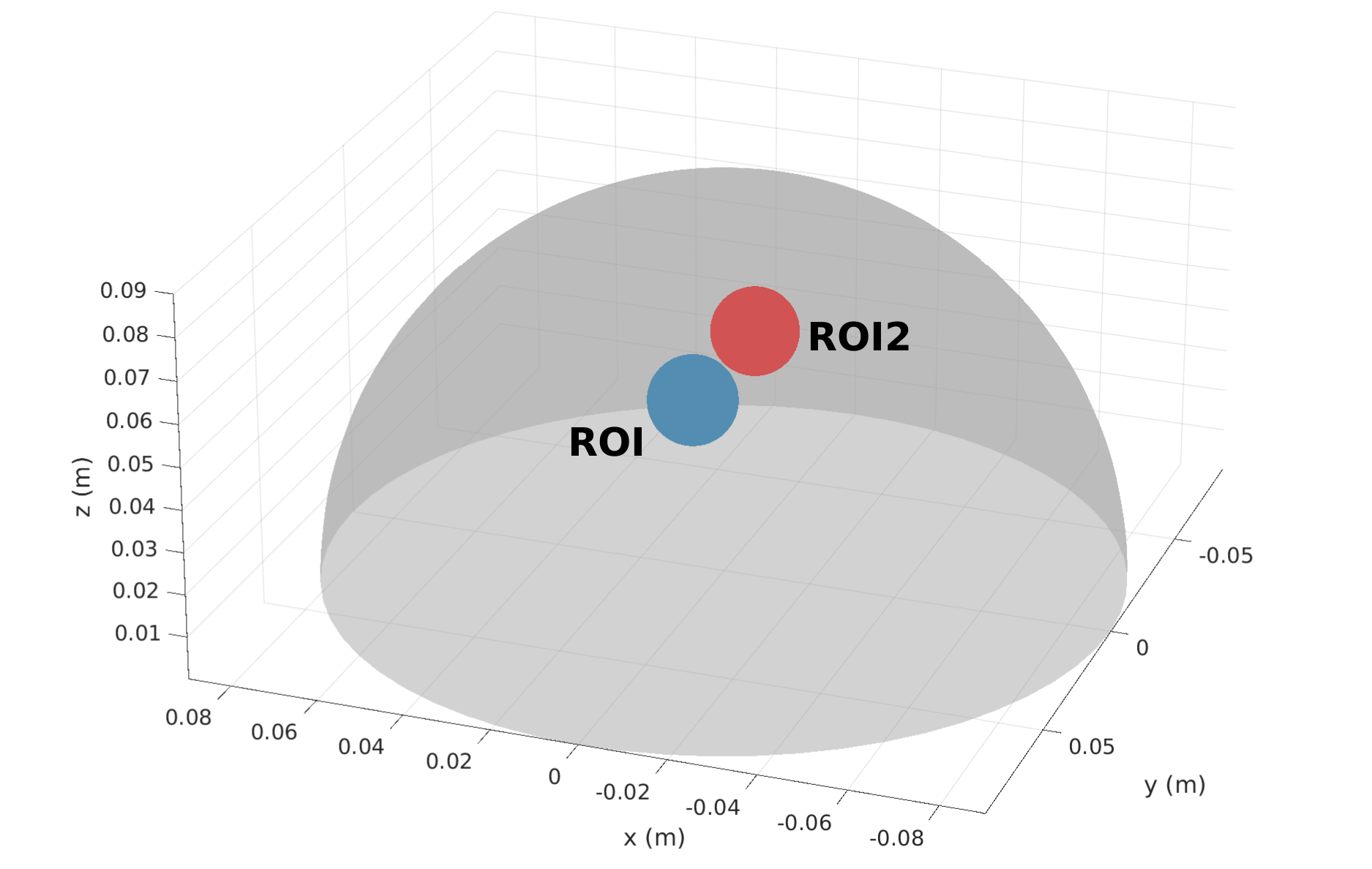}} %
\subfigure[] {\label{fig:coil}\includegraphics[width=0.9\textwidth]{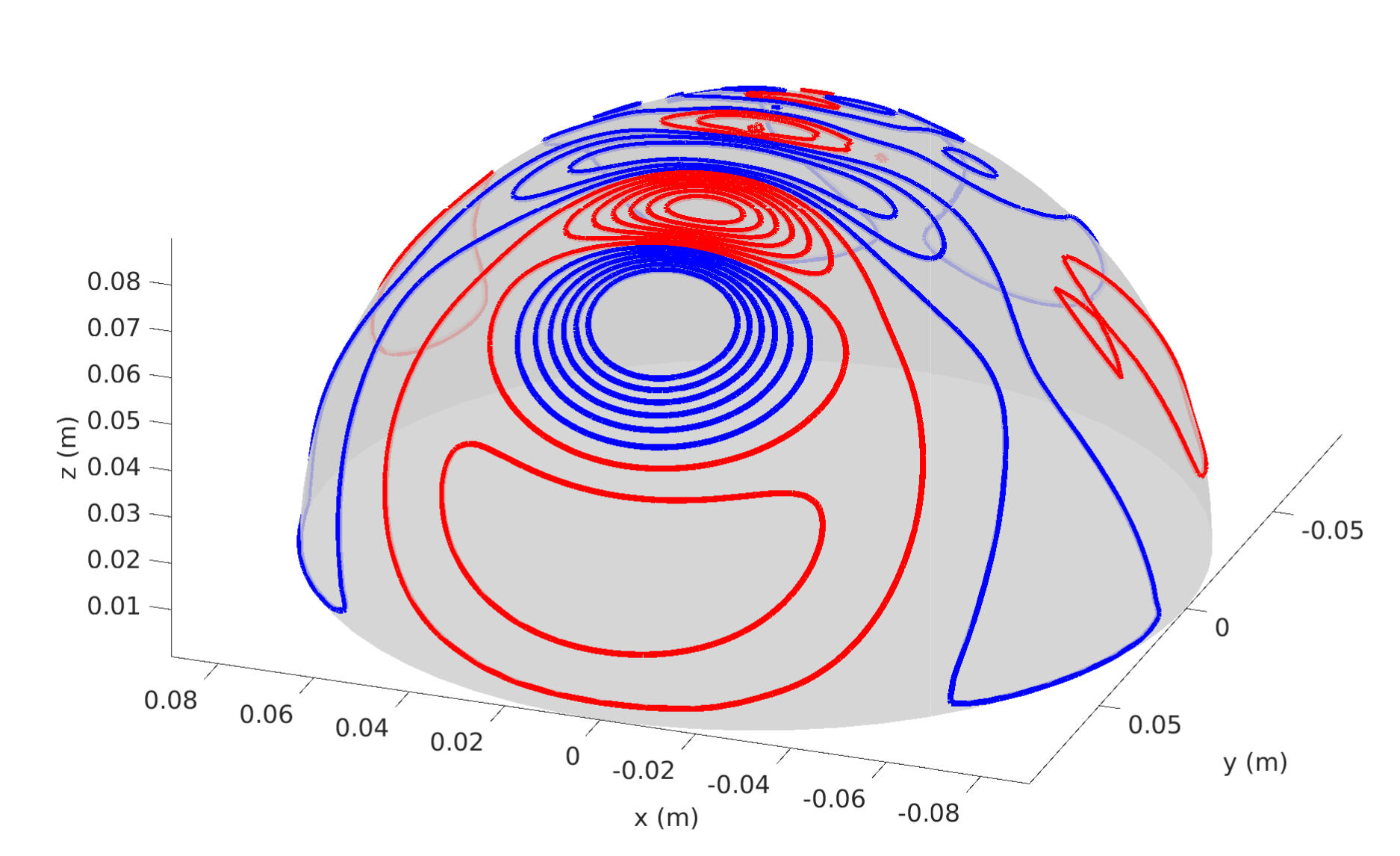}} %
\\[0pt]
\end{center}
\caption{\emph{a) Description of hemispherical surface where the optimal $\psi$ must been found along with the spherical regions of interest ROI and ROI2 where the electric field must be
maximized and minimized respectively. b) Wirepaths with 18 turns of the TMS coil solution (red wires indicate reversed current flow with respect to blue). }} \label{fig:figure1}
\end{figure}

In the following, in order to illustrate an application of the theoretical model developed in this manuscript, we are going to tackle the design of a minimum stored-energy hemispherical TMS coil of radius 9 cm, constructed to stimulate only one cerebral hemisphere. To this end, the coil must produce an E-field which is both maximum in a spherical region of interest (ROI) and  minimum in a second region (ROI2). Both volumes of interest are of 1 cm radius and  formed by 400 points, where ROI is shifted by 5 cm in the positive $z$-direction and by 2 cm in the positive $y$-direction; and ROI2 is shifted by 5 cm in the positive $z$-direction and by 2 cm in the negative $y$-direction, as shown in figure \ref{fig:rois}. By using the formalism presented in \cite{CSGPGRH} this TMS coil design problem can be posed as the following optimization problem:
\begin{equation}\label{coilTMS}
\left\{
\begin{array}{l}
\max  \left\| E_{x_1} \psi \right\| _2 \\
\min  \left\| E_{x_2} \psi \right\| _2 \\
\min   \psi^T L \psi 
\end{array}
\right.
\end{equation}
where $\psi$ is the stream function (the optimization variable), $M=400$ are the number of points in the ROI and ROI2, $N=2122$ the number of mesh nodes, $L \in \mathbb{R}^{N \times N}$ is the inductance matrix, and $E_{x_1} \in \mathbb{R}^{M \times N}$ and $E_{x_2} \in \mathbb{R}^{M \times N}$ are the $E$-field matrices in the prescribe $x$-direction.

Figure \ref{fig:coil} shows the coil solution of problem in  Eq. \ref{coilTMS} computed by using the theoretical model proposed in this manuscript (see Subsections \ref{css} and \ref{reforex}), and as expected, the wire arrangements is remarkably concentrated over the region of stimulation.

%
%

\begin{figure}[h]
\begin{center}
\subfigure[] {\label{fig:display}\includegraphics[width=0.6\textwidth]{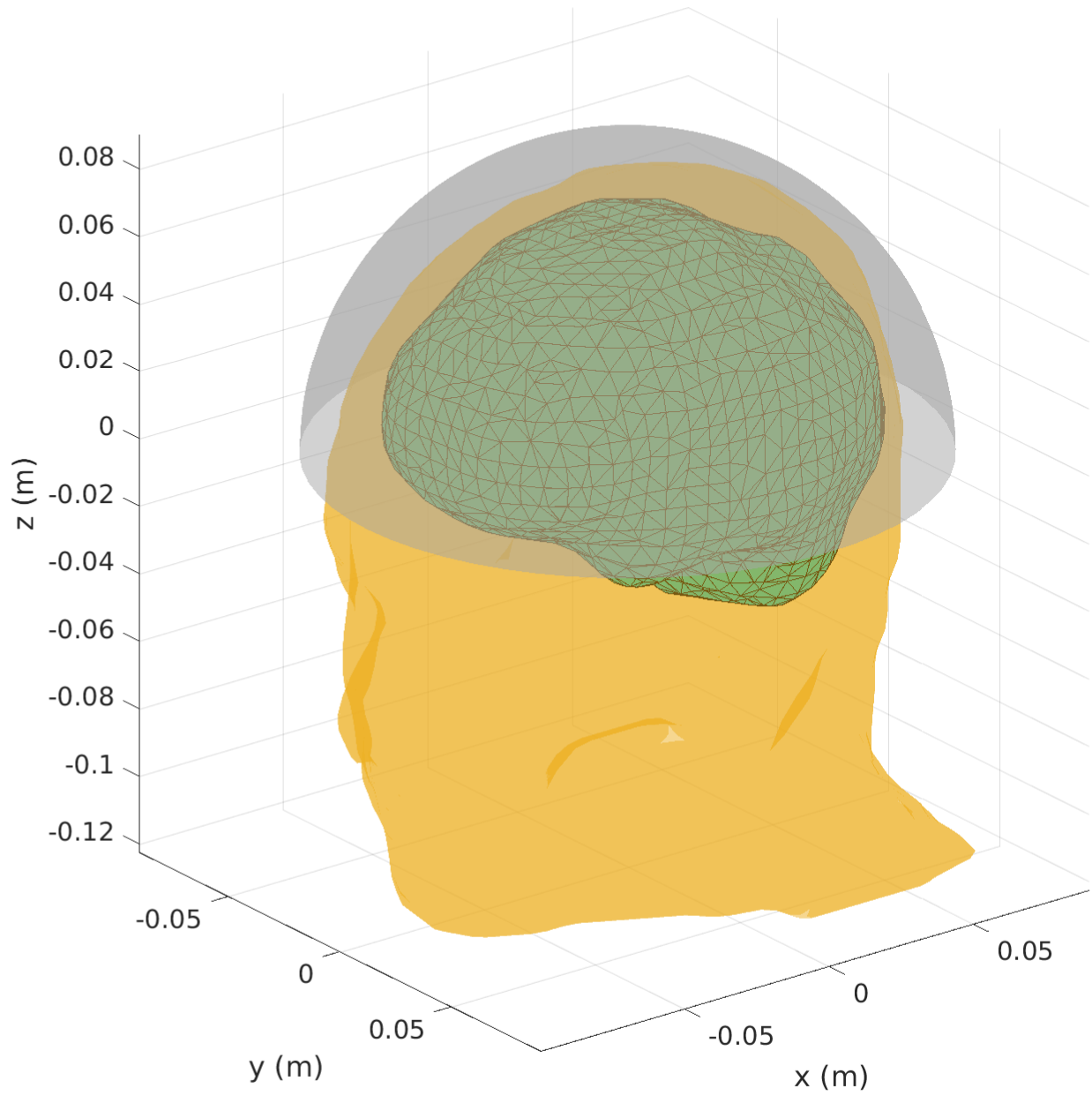}} %
\subfigure[] {\label{fig:brain}\includegraphics[width=0.5\textwidth]{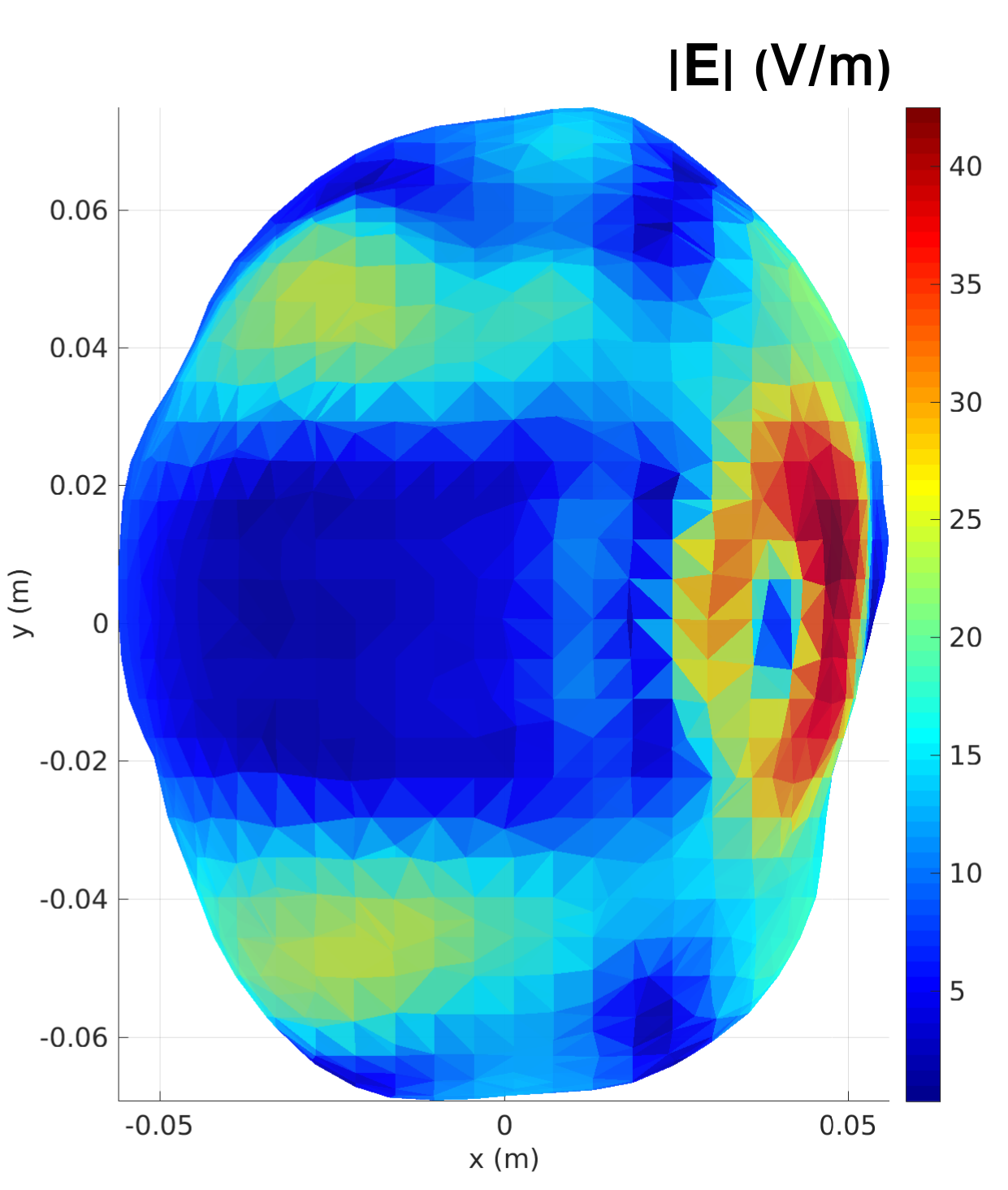}} %
\\[0pt]
\end{center}
\caption{\emph{a) Description of the two compartment scalp-brain model.  b) E-field modulus induced at the surface of the brain by the designed TMS coil. }} \label{fig:figure2}
\end{figure}

In order to evaluate the stimulation of the coil we resort to the direct BEM \cite{key:forward}, which allows calculation of the electric field induced by coils in conducting systems. In Figure \ref{fig:display} a simple human head made of two compartments, scalp and brain, used to evaluate the performance of the designed stimulator is shown. As it can be seen from Figure \ref{fig:brain}, the TMS coil fulfils the initial requirements of stimulating only one hemisphere of the brain (the one where ROI is found); whereas the electric field induced in the other cerebral hemisphere (where ROI2 can be found) is minimum.

\subsection{Reformulation of Problem \eqref{coilTMS} to turn it into a maxmin}\label{reforex}

We proceed now to reformulate the multiobjective optimization problem given in \eqref{coilTMS} in order to transform it into a maxmin problem like in \eqref{maxmint} so that we can apply the theoretical model described in Subsection \ref{css}:

\begin{equation*}
\left\{
\begin{array}{l}
\max  \left\| E_{x_1} \psi \right\| _2 \\
\min  \left\| E_{x_2} \psi \right\| _2 \\
\min   \psi^T L \psi 
\end{array}
\right.
\end{equation*}

First, by taking into consideration that raising to the square is a strictly increasing function on $[0,\infty)$, we can apply Equation \eqref{cv2} to obtain 
\begin{equation} \label{ffm1}
\left\{
\begin{array}{l}
\max  \left\| E_{x_1} \psi \right\|^2_2 \\
\min  \left\| E_{x_2} \psi \right\|^2_2 \\
\min   \psi^T L \psi\end{array}
\right.
\end{equation}
 Next, we apply Cholesky decomposition to $L$ to obtain $L=C^TC$ so we have that $\psi^TL\psi =(C\psi)^T(C\psi)=\|C\psi\|_2^2$ so we obtain
\begin{equation} \label{ffm2}
\left\{
\begin{array}{l}
\max  \left\| E_{x_1} \psi \right\|^2_2 \\
\min  \left\| E_{x_2} \psi \right\|^2_2 \\
\min   \|C\psi\|_2^2  
\end{array}
\right.
\end{equation}
 Since $C$ is an invertible square matrix, $\argmin \|C\psi\|_2^2 = \{0\}$ so the previous multiobjective optimization problem has no solution. Therefore it must be reformulated. We call then on Section \ref{moreope} to obtain:
 
 \begin{equation} \label{ffm21}
\left\{
\begin{array}{l}
\max  \left\| E_{x_1} \psi \right\|^2_2 \\
\min  \left\| E_{x_2} \psi \right\|^2_2 +   \|C\psi\|_2^2  
\end{array}
\right.
\end{equation}
which in essence is
\begin{equation} \label{ffm3}
\left\{
\begin{array}{l}
\max  \left\| E_{x_1} \psi \right\|_2 \\
\min  \left\| D \psi \right\|_2
\end{array}
\right.
\end{equation}
where $D:=\left(\begin{array}{c}
E_{x_2}\\ C\end{array}\right)$. The matrix $D$ in this specific case has null kernel. In accordance with the previous sections, Problem \eqref{ffm3} is remodeled as 
\begin{equation} \label{ffm4}
\left\{
\begin{array}{l}
\max  \left\| E_{x_1} \psi \right\|_2 \\
  \left\| D \psi \right\|_2\leq 1
\end{array}
\right.
\end{equation}

Finally, we can refer to Subsection \ref{css} to solve the latter problem.

\section{Applications to optimal geolocation}

To show another application of maxmin multiobjective problems we consider here the issue of optimal geolocation. In particular we focus on the best situation of a tourism rural inn considering several measured climate variables. Locations with low highest temperature $m_1$, radiation $m_2$ and evapotranspiration $m_3$ in summer time and high values in winter time are sites with climatic characteristics desirable for potential visitors. To solve this problem, we choose 11 locations in the Andalusian coastline and 2 in the inner, near the mountains. We have collected the data from the official {\em Andalusian government} webpage \cite{key:ifapa} evaluating the mean values of these variables on the last 5 years 2013-2019. The referred months of the study were January and July. 

\begin{table}[h!]
\begin{tabular}{l|ccc|ccc}
          & T-winter & R-winter & E-winter   & T-summer & R-summer & E-summer   \\
Sanlúcar  & 15.959      & 9.572     & 1.520     & 30.086      & 27.758    & 6.103 \\
Moguer    & 16.698      & 9.272     & 0.925     & 30.424      & 27.751    & 5.222 \\
Lepe      & 16.659      & 9.503     & 1.242     & 30.610      & 28.297    & 6.836 \\
Conil     & 16.322      & 9.940     & 1.331     & 28.913      & 26.669    & 5.596 \\
El Puerto & 16.504      & 9.767     & 1.625     & 31.052      & 28.216    & 6.829 \\
Estepona  & 16.908      & 10.194    & 1.773     & 31.233      & 27.298    & 6.246 \\
Málaga    & 17.663      & 9.968     & 1.606     & 32.358      & 27.528    & 6.378 \\
Vélez     & 18.204      & 9.819     & 1.905     & 31.912      & 26.534    & 5.911 \\
Almuñécar & 17.733      & 10.247    & 1.404     & 29.684      & 25.370    & 4.952 \\
Adra      & 17.784      & 10.198    & 1.637     & 28.929      & 26.463    & 5.143 \\
Almería   & 17.468      & 10.068    & 1.561     & 30.342      & 27.335    & 5.793 \\
Aroche    & 16.477      & 9.797     & 1.434     & 34.616      & 27.806    & 6.270 \\
Córdoba   & 14.871      & 8.952     & 1.149     & 36.375      & 28.503    & 7.615 \\
Baza      & 13.386      & 8.303     & 3.054     & 35.754      & 27.824    & 1.673 \\
Bélmez    & 13.150      & 8.216     & 1.215     & 35.272      & 28.478    & 7.400 \\
S. Yeguas & 13.656      & 9.155     & 1.247     & 33.660      & 28.727    & 7.825   
\end{tabular}
\vspace{0.5cm}
\caption{Mean values of high temperature (T) in Celsius Degrees, radiation (R) in $MJ/m^2$, and evapotranspiration (E) in mm/day, measures in January (winter time) and July (summer time) between 2013 and 2018. }
\end{table}

To find the optimal location we evaluate the site where the variables mean values are maximum in January and minimum in July. Here we have a typical multiobjective problem with two data matrices that can be formulated as follows:
\begin{equation}\label{albertproblem}
\left\{\begin{array}{l} 
\max \|Ax\|_2 \\ 
\min \|Bx\|_2\\ 
\min \|x\|_2
\end{array}\right.
\end{equation}
where $A$ and $B$ are real 16x3 matrices with the values of the three variables $(m_1,m_2,m_3)$ taking into account (highest temperature, radiation and evapotranspiration) in January and July respectively. To avoid unit effects, we standarized the variables ($\mu=0$ and $\sigma = 1$). The vector $x$ is the solution of the multiobjective problem.

This question can be reformulated as we showed in Section \ref{moreope} by the following:

\begin{equation}\label{albertproblem2}
\left\{
\begin{array}{l}
\max  \left\|Ax\right\|_2 \\
\min  \left\|Dx\right\|_2
\end{array}
\right.
\end{equation}
with matrix $D:=\left(\begin{array}{c} B\\ I_n\end{array}\right)$, where $I_n$ is the identity matrix with $n=3$. Notice that it also verifies that $\ker(D)=\{0\}$. Observe that, according to the previous sections, \eqref{albertproblem2} can be remodeled into

\begin{equation}\label{albertproblem3}
\left\{
\begin{array}{l}
\max  \left\|Ax\right\|_2 \\
  \left\|Dx\right\|_2\leq 1
\end{array}
\right.
\end{equation}

and solved accordingly.

\begin{figure}[h]
\begin{center}
\includegraphics[width=0.8\textwidth]{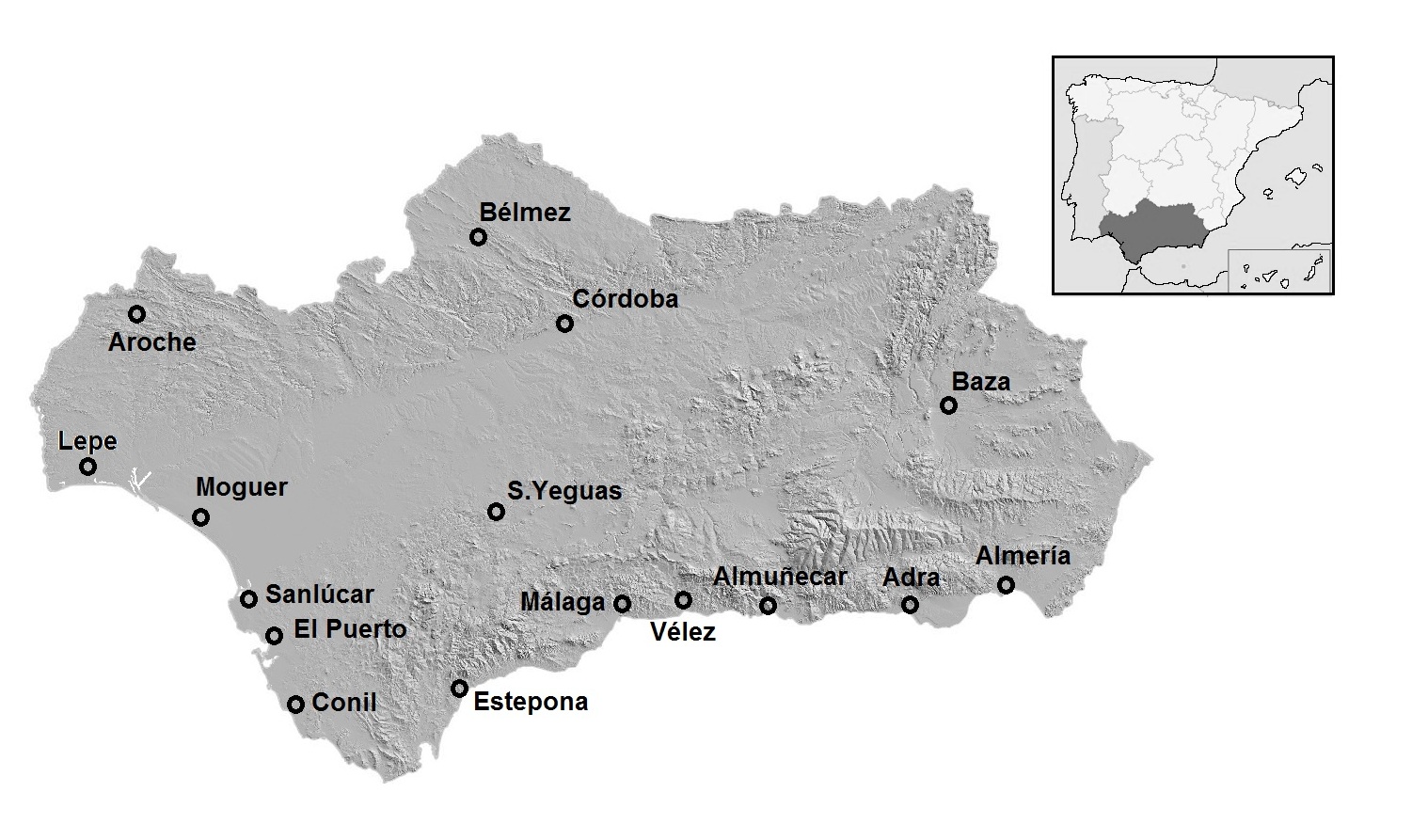}
\\[0pt]
\end{center}
\caption{\emph{Geographic distribution of the sites considered in the study. 11 places are in the coastline of the region and 5 in the inner}} \label{fig:figure3}
\end{figure}

\begin{figure}[h]
\begin{center}
\includegraphics[width=0.6\textwidth]{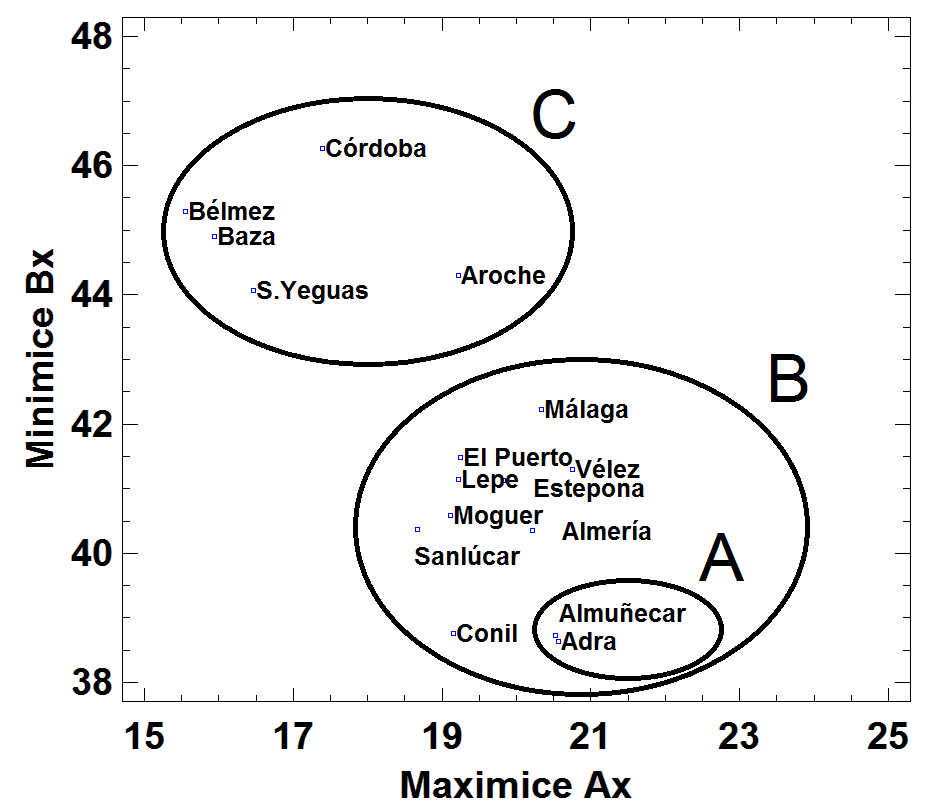}
\\[0pt]
\end{center}
\caption{\emph{Locations considering Ax and Bx axes. Group named $\it{A}$ represents the best places for the tourism rural inn, near Costa Tropical (Granada province). Sites on $\it{B}$ are also in the coastline of the region. Sites on $\it{C}$ are the worst locations considering the multiobjective problem, they are situated inside the region}} \label{fig:figure4}
\end{figure}

\begin{figure}[h]
\begin{center}
\includegraphics[width=1\textwidth]{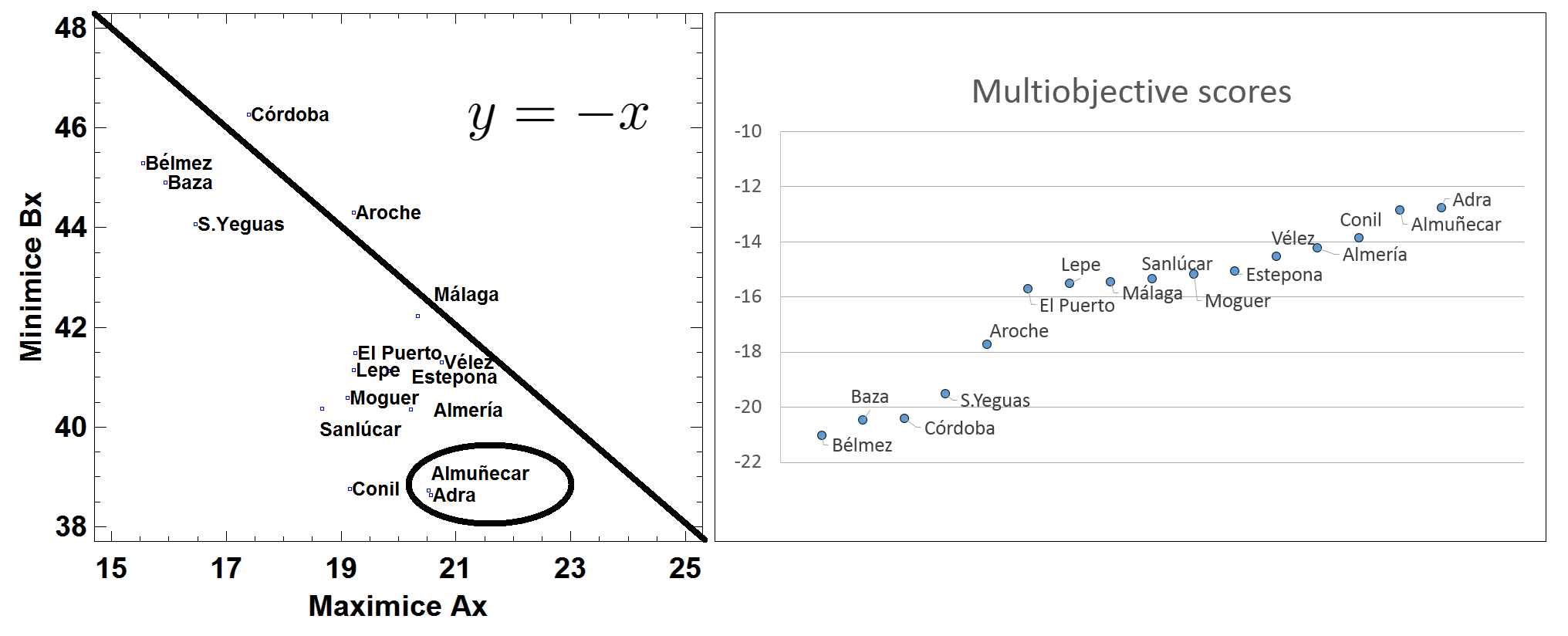}
\\[0pt]
\end{center}
\caption{\emph{a) Sites considering Ax and Bx and the function $y=-x$. The places with high values of Ax (max) and low values of Bx (min) are the best locations for the solution of the multiobjective problem (round). b) Multiobjective scores values obtained for each site projecting the point in the function $y=-x$. High values of this score indicate better places to locate the tourism rural inn.}} \label{fig:figure5}
\end{figure}

\begin{figure}[h]
\begin{center}
\includegraphics[width=1\textwidth]{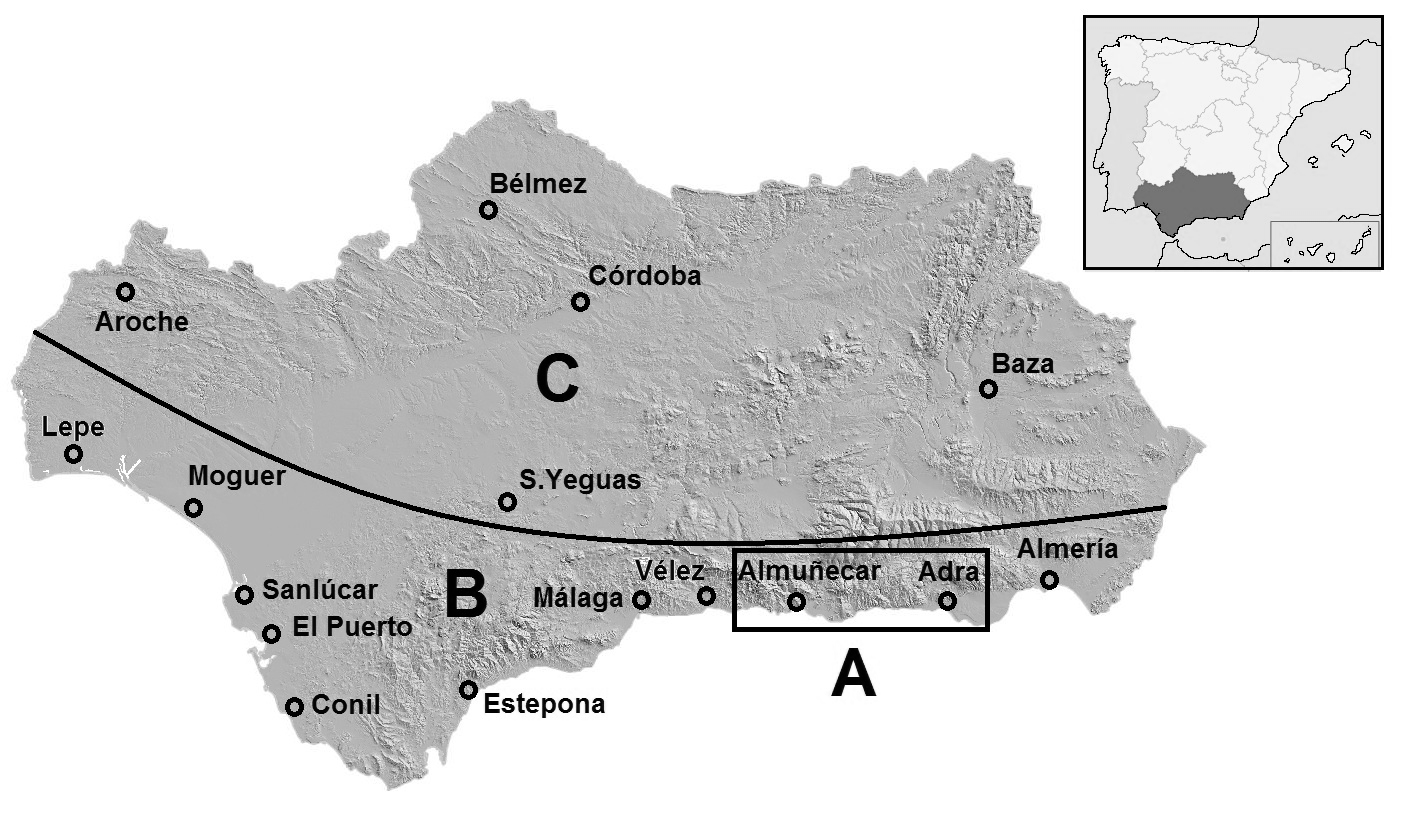}
\\[0pt]
\end{center}
\caption{\emph{Distribution of the three areas described on figure \ref{fig:figure4}. A and B areas are in the coastline and C in the inner. }} \label{fig:figure6}
\end{figure}

The solution of \eqref{albertproblem3} allow us to draw the sites with a $2D$ plot considering the $X$ axe as $Ax$ and the $Y$ axe as $Bx$. We observe that better places have high values of $Ax$ and low values of $Bx$. Hence, we can sort the sites in order to achieve the objectives in a similar way as factorial analysis works (two factors, the maximum and the minimum, instead of $m$ variables).

\section{Algorithms developed in this work}

In this section we show the algorithms written in MATLAB to solve the real problems presented. First we include the solution (generalized supporting vectors) of the problem presented in \cite{CSGPMPSM} which is (this code appears in \cite{CSGPMPSA}):
\begin{align*}
    \max_{\|x\|_2 = 1} \sum_{i=1}^k \| A_i x\|_2^2 = \lambda_{\max} \left(\sum_{i=1}^k A_i^T A_i\right) \\
    \arg\max_{\|x\|_2 = 1} \sum_{i=1}^k \| A_i x\|_2^2  = V\left(\lambda_{\max}\left(\sum_{i=1}^k A_i^T A_i\right)\right) \cap \E_{\ell_2^n}
\end{align*}
where $\lambda_{\max}$ denotes the greatest eigenvalue and $V$ denotes the associated eigenvector space.
\begin{lstlisting}
    function [lambda_max, x] = sol_1(M)
    %%%%%
    %%%%% INPUT:
    %%%%%
    %%%%% M = {A_1 A_2 ... A_k} a list with the matrices
    %%%%%
    %%%%%%%%%%%%%%%%%%%%%%%%%%%%%%%%%
    %%%%%
    %%%%% OUTPUTS:
    %%%%%
    %%%%% lambda_max - maximum eigenvalue
    %%%%% x - eigenvector associated to lambda_max
    %%%%%%%%%%%%%%%%%%%%%%%%%%%%%%%%%
    %%%%
    k = length(M);
    [nrows ncols] = size(M{1});
    a1 = M{1}(:,1); %% First column of M{1}
    a2 = M{1}(:,2); %% Second column of M{1}
    if (k==1) & (ncols==2) & (abs(norm(a1)-norm(a2))<1e-12) 
        %%% In this particular case, a single matrix with two columns with 
        %%% the same norm is considered. A tolerance of 1e-12 is needed in 
        %%% order to compare these norms.
        %%% The maximum lambda_max and the supporting vectors are computed 
        %%% directly 
        lambda_max = norm(a1)^2 + abs(a1'*a2);
        if floor(a1'*a2)==0
            %%% The columns of this matrix form a basis of supporting vectors
            x = eye(2); 
        elseif a1'*a2>0
            %%% The columns of this matrix form a basis of supporting vectors
            x = [sqrt(2)/2 sqrt(2)/2; -sqrt(2)/2 -sqrt(2)/2]'; 
        elseif a1'*a2<0
            %%% The columns of this matrix form a basis of unit supporting vectors
            x = [-sqrt(2)/2 sqrt(2)/2; sqrt(2)/2 -sqrt(2)/2]';
        end
    else
        %%% This is the general case
        suma = zeros(ncols);
        for i=1:k
            suma = suma + M{i}'*M{i};
        end
        %%%
        [V,D] = eig(suma);           %%% Computing the eigensystem
        lambda_max = max(diag(D)); %%% This is the maximum eigenvalue
        N = size(D,1);
        %%% Now we find the indices where the elements of the diagonal of the
        %%% matrix D are equal (with a tolerance of 1e-12) to lambda_max  
        ind_lambda_max = find(abs(diag(D)-lambda_max*ones(N,1))<1e-12);
        x = V(:, ind_lambda_max); %%% %%% The columns of this matrix form a
        %%% basis of unit supporting vectors associated to the maximum eigenvalue 
    end
end
\end{lstlisting}
As we pointed out in Theorem \ref{maxminchar}, the solution of the problem 
$$\left\{\begin{array}{l} \max \|Ax\| \\ \|Bx\|\leq 1 \end{array}\right.$$
exists if and only if $\ker(B)\subseteq \ker(A)$. Here is a simple code to check this.
\begin{lstlisting}
    function p=existence_sol(A,B)
    %%%%
    %%%% This function checks the existence of the solution of the
    %%%% problem
    %%%% 
    %%%% max ||Ax||
    %%%% ||Bx||<=1
    %%%%
    %%%%%%%%%%%%%%%%%%%%%%%%%%%%%%%
    %%%%
    %%%% INPUT:
    %%%%
    %%%% A, B - the matrices involved in the problem
    %%%%
    %%%%%%%%%%%%%%%%%%%%%%%%%%%%%%%%
    %%%%
    %%%% OUTPUT:
    %%%%
    %%%% p - true if the problem has solution or false on the contrary
    %%%% 
    %%%%%%%%%%%%%%%%%%%%%%%%%%%%%%%%
    KerB = null(B);
    dimKerB = size(KerB,2);
    KerA = null(A);
    dimKerA = size(KerA,2);
    if (dimKerB<=dimKerA) & (rank([KerB KerA])==dimKerA)
        p = true;
    else
        p = false;
    end
end
\end{lstlisting}
Now we present the code to solve the first case of the previous maxmin problem, that is, the case where $\ker(B)=\{0\}$. We refer the reader to Subsection \ref{css} upon which this code is based.

\begin{lstlisting}
function x = case_1(A, B)
    %%%%
    %%%% This function computes the solution of the problem 
    %%%% 
    %%%% max ||Ax||_2
    %%%% ||Bx||_2<=1
    %%%%
    %%%% in the case KerB={0}.
    %%%%%%%%%%%%%%%%%%%%%%%%%%%%%%%
    %%%%
    %%%% INPUT:
    %%%%
    %%%% A, B - the matrices involved in the problem
    %%%%
    %%%%%%%%%%%%%%%%%%%%%%%%%%%%%%%%
    %%%%
    %%%% OUTPUT:
    %%%%
    %%%% x - basis of unit  eigenvectors associated to lambda_max
    %%%% 
    %%%%%%%%%%%%%%%%%%%%%%%%%%%%%%%%
    %%%%
    KerB = null(B);
    dimKerB = size(KerB,2);
    if (dimKerB ~= 0)
        display('KerB~={0}')
        x=[];
    else % KerB={0}
            M = A*pinv(B);              % M = A*B^+
                                        % B^+ is the pseudoinverse matrix
            [lambda_max, y] = sol_1({M});
            [nrows_y ncols_y] = size(y);
            r_B = rank(B);
            counter = 0;
            for i=1:ncols_y
                r = rank([B y(:,i)]);
                if (abs(r_B - r)<1e-12)     % Here we check if rank(B) = rank ([B y0]). 
                                % A tolerance of 1e-12 is needed in 
                                % order to compare these two ranks.
                   counter = counter +1;
                   y0(:,counter) = y(:,i);
                end
            end
            x = pinv(B)*y0;               % This is a basis of solutions of our problem
end
\end{lstlisting}

Next, we can compute the global solution of the maxmin problem by means of the following code. Again, we refer the reader to Subsection \ref{css} upon which this code is based.
\begin{lstlisting}
    function x = sol_2(A, B)
    %%%%
    %%%% This function computes the solution of the problem 
    %%%% 
    %%%% max ||Ax||_2
    %%%% ||Bx||_2<=1
    %%%%
    %%%%%%%%%%%%%%%%%%%%%%%%%%%%%%%
    %%%%
    %%%% INPUT:
    %%%%
    %%%% A, B - the matrices involved in the problem
    %%%%
    %%%%%%%%%%%%%%%%%%%%%%%%%%%%%%%%
    %%%%
    %%%% OUTPUT:
    %%%%
    %%%% x - Supporting vector which is the solution of the problem
    %%%% 
    %%%%%%%%%%%%%%%%%%%%%%%%%%%%%%%%
    %%%%
    p=existence_sol(A,B);
    if p==true                        
        n = size(B,2);
        KerB = null(B);
        dimKerB = size(KerB,2);
        if (dimKerB == 0)               % KerB = {0} This is the case 1
            x = case_1(A,B);             % x is the solution of our problem
        else % KerB~={0}
            [Br indices] = colsindep(B); %%% First we extract the 
                                            %%% independent columns in B
            Ar = A(indices);        %%% We extract the same columns of A
                  %%% Now, Ker(Br)={0} so this is the case 1 treated above:
            xr = case_1(Ar,Br);
            [nrows_xr,ncols_xr] = size(xr);
            %%% Now we compute the matrix solutions x of the problem
            counter = 0;
            for j = 1:ncols_xr
                for i=1:n
                    if ismember(i,indices)==1 %%% i is an index of the ones 
                                          %%% defined above
                        counter = counter + 1;
                        x(i,j) = xr(counter,j);
                    else
                        x(i,j) = 0;
                    end
                end
            end
        end
        
    else
        display('This problem has no solution');
        x=[];
    end
end
\end{lstlisting}
Notice that we use the \verb case_1 \mbox{ }function described above and a new function named \verb colsindep. We include the code to implement this new function below.
\begin{lstlisting}
    function [Dcolsind, indices]=colsindep(D)
    %%%%
    %%%% This function extracts r = rank(D) independent columns of the
    %%%% matrix D and the indices of the columns in D which are independent 
    %%%%
    %%%%%%%%%%%%%%%%%%%%%%%%%%%%%%%
    %%%%
    %%%% INPUT:
    %%%%
    %%%% D - a matrix with rank r
    %%%%
    %%%%%%%%%%%%%%%%%%%%%%%%%%%%%%%%
    %%%%
    %%%% OUTPUT:
    %%%%
    %%%% Dcolsind - r independent columns in D
    %%%% indices - the indices of independent columns extracted from D
    %%%%%%%%%%%%%%%%%%%%%%%%%%%%%%%%
    r=rank(D);            %%% Compute the rank
    [Q R p]=qr(D,0);      %%% p is a permutation vector such that A(:,p)=Q*R
    indices=sort(p(1:r)); %%% The first r elements in p are the indices of the
                          %%% columns linearly independent in D
    Dcolsind=D(:,indices);%%% Extract these columns
end
\end{lstlisting}

Here we include the code to compute the solution of the TMS coil problem \eqref{ffm4}:

\begin{equation*}
\left\{
\begin{array}{l}
\max  \left\| E_{x_1} \psi \right\|_2 \\
  \left\| D \psi \right\|_2\leq 1
\end{array}
\right.
\end{equation*}
with the matrix $D:=\left(\begin{array}{c} E_{x_2}\\ C\end{array}\right)$, where $C$ is the Cholesky matrix of $L$, and in this case it verifies that $\ker(D)=\{0\}$. Recall that \eqref{ffm4} comes from \eqref{coilTMS}:
\begin{equation*}
\left\{
\begin{array}{l}
\max  \left\| E_{x_1} \psi \right\| _2 \\
\min  \left\| E_{x_2} \psi \right\| _2 \\
\min   \psi^T L \psi  
\end{array}
\right.
\end{equation*}

\begin{lstlisting}
function psi = sol2_psi(Ex1, Ex2, L)

    C = chol(L);                    % Cholesky's decomposition of matrix L = C' * C

    A = Ex1;
    B = [Ex2;C];                 

    psi = case_1(A,B);           % We apply the algorithm to obtain the solutions
end
\end{lstlisting}

Finally, we provide the code to compute the solution of the optimal geolocation problem \eqref{albertproblem3}:

\begin{equation*}
\left\{
\begin{array}{l}
\max  \left\|Ax\right\|_2 \\
 \left\|Dx\right\|_2\leq 1
\end{array}
\right.
\end{equation*}
with matrix $D:=\left(\begin{array}{c} B\\ I_{3}\end{array}\right)$. Notice that it also verifies that $\ker(D)=\{0\}$ and $A$ and $B$ are composed by standardized variables. Recall that \eqref{albertproblem3} comes from \eqref{albertproblem}:

$$\left\{\begin{array}{l}
\max \|Ax\|_2 \\ 
\min \|Bx\|_2\\ 
\min \|x\|_2
\end{array}\right.$$

\begin{lstlisting}
function x = sol_2_geoloc(A, B)
    
    [rows,cols] = size(A);
    D = [B; eye(size(cols))];                    
    
    x = case_1(A,D);           % We apply the algorithm to obtain the solutions
end
\end{lstlisting}


\begin{thebibliography}{10}

\bibitem{Bo}
Frederic Bohnenblust.
\newblock A characterization of complex {H}ilbert spaces.
\newblock {\em Portugal. Math.}, 3:103--109, 1942.

\bibitem{CSGPGRH}
Clemente Cobos~S\'{a}nchez, Francisco~Javier Garcia-Pacheco, Jose~Maria
  Guerrero~Rodriguez, and Justin~Robert Hill.
\newblock An inverse boundary element method computational framework for
  designing optimal {TMS} coils.
\newblock {\em Eng. Anal. Bound. Elem.}, 88:156--169, 2018.

\bibitem{CSGPMPSM}
Clemente Cobos-S\'{a}nchez, Francisco~Javier Garc\'{\i}a-Pacheco, Soledad
  Moreno-Pulido, and Sol S\'{a}ez-Mart\'{\i}nez.
\newblock Supporting vectors of continuous linear operators.
\newblock {\em Ann. Funct. Anal.}, 8(4):520--530, 2017.

\bibitem{key:wassermann}
Ulf~Ziemann Eric~Wassermann, Charles~Epstein.
\newblock {\em Oxford Handbook of Transcranial Stimulation (Oxford Handbooks)}.
\newblock Oxford University Press, New York, 1 edition, 2008.

\bibitem{CSGPMPSA}
Francisco~Javier Garcia-Pacheco, Clemente Cobos-Sanchez, Soledad Moreno-Pulido,
  and Alberto Sanchez-Alzola.
\newblock Exact solutions to {$\max_{\Vert x\Vert=1}\sum^\infty_{i=1}\Vert
  T_i(x)\Vert^2$} with applications to {P}hysics, {B}ioengineering and
  {S}tatistics.
\newblock {\em Commun. Nonlinear Sci. Numer. Simul.}, 82:105054, 2020.

\bibitem{GPNG}
Francisco~Javier García-Pacheco and Enrique Naranjo-Guerra.
\newblock Supporting vectors of continuous linear projections.
\newblock {\em International Journal of Functional Analysis, Operator Theory
  and Applications}, 9(3):85--95, 2017.

\bibitem{H}
Na~Huang and Chang-Feng Ma.
\newblock Modified conjugate gradient method for obtaining the minimum-norm
  solution of the generalized coupled {S}ylvester-conjugate matrix equations.
\newblock {\em Appl. Math. Model.}, 40(2):1260--1275, 2016.

\bibitem{Ka}
S.~Kakutani.
\newblock Some characterizations of {E}uclidean space.
\newblock {\em Jap. J. Math.}, 16:93--97, 1939.

\bibitem{key:forward}
Clemente~Cobos Sanchez, Richard~W. Bowtell, Henry Power, Paul Glover, Liviu
  Marin, Adib~A. Becker, and Arthur Jones.
\newblock Forward electric field calculation using {BEM} for time-varying
  magnetic field gradients and motion in strong static fields.
\newblock {\em Eng. Anal. Bound. Elem.}, 33(8-9):1074--1088, 2009.

\bibitem{key:ifapa}
\url{https://www.juntadeandalucia.es/agriculturaypesca/ifapa/ria/servlet/FrontController}.

\bibitem{Y}
Belkourchia Yassin, Azrar Lahcen, and Es-Sadek~Mohamed Zeriab.
\newblock Hybrid optimization procedure applied to optimal location finding for
  piezoelectric actuators and sensors for active vibration control.
\newblock {\em Applied Mathematical Modelling}, 62:701 -- 716, 2018.

\end{thebibliography}
\end{document}